\documentclass[reqno]{amsart}

\usepackage[dvipdfm,truedimen,left=25truemm,right=25truemm,top=30truemm,bottom=30truemm]{geometry}
\usepackage{amsfonts}
\usepackage{amsmath,amsthm}
\usepackage{amsfonts,amssymb}
\usepackage{mathrsfs}
\usepackage{enumerate}
\usepackage{bbm}
\usepackage{array}

\allowdisplaybreaks

\newtheorem{thm}{Theorem}[section]

\newtheorem{lem}[thm]{Lemma}
\newtheorem{prop}[thm]{Proposition}

\theoremstyle{remark}
\newtheorem{rem}[thm]{Remark}

\numberwithin{equation}{section}

\begin{document}

\title[] {A positive product formula of integral kernels of $k$-Hankel transforms}

\author[]{Wentao Teng} 

\address{Graduate School of Mathematical Sciences, The University of Tokyo, 3-8-1 Komaba Meguro-ku Tokyo 153-8914, Japan}
\email{wentaoteng6@sina.com.}

\begin{abstract} 
The $k$-Hankel transform $F_{k,1}$ (or the $(k,1)$-generalized Fourier transform) is the Dunkl analogue of the unitary inversion operator in the minimal representation of a conformal group initiated by T. Kobayashi and G. Mano.
It is one of the two most significant cases in $(k,a)$-generalized Fourier transforms. 
We will establish a positive radial product formula for the integral kernels of $F_{k,1}$.
Such a product formula is equivalent to a representation of the generalized spherical mean operator in terms of the probability measure $\sigma_{x,t}^{k,1}(\xi)$.
We will then study the representing measure $\sigma_{x,t}^{k,1}(\xi)$ and analyze the support of this measure, and derive a weak Huygens's principle for the deformed wave equation in $(k,1)$-generalized Fourier analysis.
\end{abstract}




\maketitle
\input amssym.def

\section{Introduction}

Let $\triangle$ be the ordinary Euclidean Laplacian on $\mathbb R^N$. For the classical Fourier transform 
$$F(f)(y)=(2 \pi)^{-N / 2} \int_{\mathbb{R}^N} f(x) e^{-i\langle x, y\rangle} d x,$$ 
R. Howe \cite{H} gave the following spectral definition of $F$ using the harmonic oscillator $\mathbf H=:(\triangle-\left|x\right|^2)/2$ and its eigenfunctions forming the basis in ${{L}^{2}}\left( {{\mathbb{R}}^{N}}\right)$:
\begin{align}\label{Howe}F:=e^{i\pi N/4}\exp\left(\frac{\mathrm{\pi i}}2\mathbf H\right).\end{align}

Various kinds of generalizations of the classical Fourier transform have drawn considerable attention
 during the last 30 years. 
One example was the Dunkl transform $F_k$, which was given in \cite{Du2} and defined with the help of a root system $R\subset\mathbb{R}^N$, a finite reflection group $G$, and a $G$-invariant multiplicity function $k:R\rightarrow \mathbb C$. 
The finite reflection group $G$ plays the role of the orthogonal group $O(N)$ in classical Fourier analysis. When $k\equiv0$, the Dunkl transform $F_k$ reduces to the classical Fourier transform $F$. The differential-difference operator $\triangle_k$, called Dunkl Laplacian, takes the place of classical Laplacian in classical analysis. It plays an important role in Dunkl analysis. If $k\equiv0$, we have $\triangle_k=\triangle$.

Motivated by the definition \eqref{Howe} of the classical Fourier transform on
${{L}^{2}}\left( {{\mathbb{R}}^{N}}\right)$ by Howe, S. Ben Sa\"id, T. Kobayashi and B. \O rsted
\cite{BSKO} gave a further far-reaching generalization of Dunkl
transform $F_k$ by introducing a parameter $a>0$ arisen from the
``interpolation'' of the two $\mathfrak{sl}(2,\mathbb R)$ actions
on the Weil representation of the metaplectic group $Mp(N,\mathbb
R)$ and the minimal unitary representation of the conformal group
$O(N+1,2)$.
They defined an $a$-deformed Dunkl harmonic oscillator  
$\triangle_{k,a}:=\left|x\right|^{2-a}\triangle_k-\left|x\right|^a$ on the Hilbert space ${{L}^{2}}\left(
{{\mathbb{R}}^{N}},{\vartheta_{k,a}}\left( x \right)dx \right)$
 with a dense domain $W_{k,a}(\mathbb R^N)$ (see \cite[Section 3]{BSKO}), and the $(k,a)$-generalized Fourier transform
$$F_{k,a}=e^{i\pi(\frac{2\left\langle k\right\rangle+N+a-2}{2a})}\exp\left(\frac{\pi i}{2a}\triangle_{k,a}\right).$$
Here ${\vartheta_{k,a}}\left( x \right)={{\left| x
\right|}^{a-2}}\vartheta_k(x)$ and $\vartheta_k\left(x\right)=\prod_{\alpha\in R}\vert\left\langle\alpha,\;x\right\rangle\vert^{k(\alpha)}$.
The $(k,a)$-generalized Fourier transform $F_{k,a}$ is a unitary operator on ${{L}^{2}}\left(
{{\mathbb{R}}^{N}},{\vartheta_{k,a}}\left( x \right)dx \right)$ with the norm 
$${\left\|f\right\|}_{2,\vartheta_{k,a}(x)dx}=\left(\int_{\mathbb{R}^N}\vert f(x)\vert^2\vartheta_{k,a}(x)dx\right)^{1/2}.$$
When $a=2$, it reduces to the Dunkl transform $F_k$.
The $(k,a)$-generalized Fourier
transform has the following integral representation on ${{L}^{2}}\left(
{{\mathbb{R}}^{N}},{\vartheta_{k,a}}\left( x \right)dx \right)$ (see \cite[(5.8)]{BSKO}) by Schwartz kernel theorem
$$F_{k,a}f\left(\xi\right)=c_{k,a}\int_{\mathbb{R}^N}f\left(y\right)B_{k,a}\left(\xi,y\right)\vartheta_{k,a}\left(y\right)dy,\;\;\xi\in\mathbb{R}^N,$$
where $c_{k,a}$ is a constant.
The two special cases for $a=2$ (the Dunkl case) and $a=1$ are of particular interest, since $(k,a)$-generalized Fourier analysis is the ``interpolation'' of the two special cases, and they bring up much richer structures in the generalization of Fourier analysis and more applications to quantum many body systems, random matrices, and many other problems mathematical physics.
While the case of $a=2$, known as Dunkl theory, has been intensively  studied during the past thirty years, the study of the case $a=1$, which was initiated by T. Kobayashi and G. Mano \cite{KM1} as a projection of the Fourier analysis on an isotropic cone (see \cite{KM2}), is still at its infancy and it inherits the conic structure.
We will focus on the special case of $a=1$, in which the generalized Fourier transform is the Dunkl generalization of unitary inversion operator in the minimal unitary representation of the conformal group
$O(N+1,2)$.

For $a={\textstyle\frac2{n}},\;n\in{\mathbb{N}}_+$ and under certain conditions for $a$, $N$ and the multiplicity function $k$, the integral kernel $B_{k,a}(x,y)$ is uniformly bounded and we have the inversion formulas of the $(k,a)$-generalized Fourier transforms. In such cases we can define the $(k,a)$-generalized translations $\tau_x$ via the generalized Fourier transform. 
And for the two particular cases $a=1$ and $a=2$ (the Dunkl case) assuming that $2\left\langle k\right\rangle+N+a-3>0$, we have the formulas of the $(k,a)$-generalized translation operator for radial functions. The radial formula for $a=2$ was found by R\"osler \cite{R2} and for $a=1$ it was found by S. Ben Sa\"id and L. Deleaval \cite{BD}. 
For both the two cases the generalized translation operators are positive on radial functions but not necessarily positivity-preserving on non-radial functions.

In \cite {BNS}, M. A. Boubatra, S. Negzaoui and M. Sifi established the following product formulas for the one dimensional case $N=1$ of the integral kernels of the $(k,a)$-generalized Fourier transform for $a={\textstyle\frac2n},\;n\in{\mathbb{N}}_+$, $\left\langle k\right\rangle>{\textstyle\frac12}-{\textstyle\frac a4}$,
\begin{align} \label{product formula}B_{k,a}(x,z)B_{k,a}(y,z)=\int_{\mathbb{R}}B_{k,a}(\xi,z)d\nu_{x,y}^{k,a}(\xi),\;z\in\mathbb{R},\end{align}
where the measures $d\nu_{x,y}^{k,a}$ are uniformly bounded signed Borel measures on $\mathbb{R}$.
They are equivalent to the following integral representations of the $(k,a)$-generalized translations
\begin{align} \label{tauN=1}\tau_xf(y)=\int_{\mathbb{R}}f(\xi)d\nu_{x,y}^{k,a}(\xi),\;a={\textstyle\frac2{n}},\;n\in{\mathbb{N}}.\end{align}
The Borel measures $d\nu_{x,y}^{k,a}$, however, are not positive in
contrast to a hypergroup convolution.
They satisfy $\nu_{x, y}^{k, a}(\mathbb{R})=1$ and that $\operatorname{supp}\left(\nu_{x,y}^{k,a}\right)(\mathbb{R})\subset\{z\in\mathbb{R}/\left.\vert\vert x\right|^{\textstyle\frac1n}-\vert y{\vert^{\textstyle\frac1n}\vert}<\left|z\right|^{\textstyle\frac1n}<\vert x\vert^{\textstyle\frac1n}+\vert y\vert^{\textstyle\frac1n}\}$.
And in particular, the special case of the product formula for $a=1$ was given in \cite{Be} by S. Ben Sa\"id.
It could be conjectured that the above product formulas hold true for arbitrary dimensions when $a={\textstyle\frac2n},\;n\in{\mathbb{N}}_+$, which has been a long open conjecture for $a=2$ (see \cite{R1}).

Consider the spherical mean operator $f\mapsto M_f,\;f\in C(\mathbb{R}^N)$ of the $(k,a)$-generalized translation $\tau_x$ for $a=2$ (the Dunkl case) and $1$,
where 
$$M_f(x,t):=\frac1{d_{k,a}}\int_{\mathbb S^{N-1}}\tau_xf(ty)\vartheta_{k,a}(y)d\sigma(y),\;d_{k,a}=\int_{\mathbb S^{N-1}}\vartheta_{k,a}(y)d\sigma(y),\;\quad(x\in\mathbb{R}^N,\,t\geq0).$$ 
The generalized translation $\tau_x$ is usually not positive, but it was shown in \cite{R2} that the spherical mean operator for $a=2$ is positivity-preserving and it is uniquely represented by a compactly supported probability measure.
We consider the particular case for $a=1$ in this paper.

\begin{thm}\label{positivity}
1). Assume $2\left\langle k\right\rangle+N-2>0$. For $a=1$, the spherical mean operator $f\mapsto M_f$ is positivity-preserving on $C_0(\mathbb{R}^N)$, i.e., if $f\in C_0(\mathbb{R}^N)$ and $f\geq 0$ on $\mathbb R^N$, then $M_f\geq 0$ on $\mathbb R^N\times \mathbb R^+$.\\
2). Under the conditions in 1), the operator $f\mapsto M_f$  is represented by a unique probability measure $\sigma_{x,t}^{k,1}\in M^1(\mathbb{R}^N)$, i.e.,
\begin{align}\label{Mfxt} M_f(x,t)=\int_{\mathbb{R}^N}f\,d\sigma_{x,t}^{k,1}\quad\text{for all }\,f\in C_0(\mathbb{R}^N),\end{align}
where $M^1(\mathbb{R}^N)$ stands for the space of Borel probability measures.
The measure $\sigma_{x,t}^{k,1}$ has a compact support and it satisfies
$$\operatorname{supp}\sigma_{x,t}^{k,1}\subseteq\left\{\xi\in\mathbb{R}^N:\sqrt{\vert\xi\vert}\geq\vert\sqrt{\vert x\vert}-\sqrt t\vert\right\}\cap\left(\bigcup_{g\in G}\{\xi\in\mathbb{R}^N:d(\xi,gx)\leq\sqrt t\}\right),$$
where $d\left(x,y\right)=\sqrt{\left|x\right|+\left|y\right|-\sqrt{2\left(\left|x\right|\left|y\right|+\left\langle x,y\right\rangle\right)}}$. Moreover,
\begin{align} \label{transformation}
\sigma_{g x, t}^{k,1}(A)=\sigma_{x, t}^{k,1}\left(g^{-1}(A)\right) \quad \text { and } \quad \sigma_{r x, r t}^{k,1}(A)=\sigma_{x, t}^{k,1}\left(r^{-1} A\right)
\end{align}
for all $g \in G, r>0$, and all Borel sets $A \in \mathcal{B}\left(\mathbb{R}^N\right)$. \\
3). 
The formula \eqref{Mfxt} can be extended to all bounded continuous functions $ C_b(\mathbb{R}^N)$ since the representing measure $\sigma_{x, t}^{k,1}$ is compactly supported according to 2). And  
\eqref{Mfxt} is equivalent to the following product formula
\begin{align}\label{product formula2}B_{k,1}(x,z)j_{2\left\langle k\right\rangle+N-2}\left(2\sqrt{\left|tz\right|}\right)=\int_{\mathbb{R}^N}B_{k,1}(\xi,z)\,d\sigma_{x,t}^{k,1}(\xi)\; \;\;\;\;\mathrm{for}\;\mathrm{all}\;z\in\mathbb{R}^N,\end{align}
where $j_{\lambda}(x)$ is the normalized Bessel function.
\end{thm}

\begin{rem}

Such a product formula \eqref{product formula2} is a ``radial'' form of the product formula  \eqref{product formula} for $a=1$ and arbitrary dimensions, and can also be regarded as a stronger version of the positive radial formula (ref. \eqref{radial}) of $(k,1)$-generalized translations, since $j_{2\left\langle k\right\rangle+N-2}(x)$ is the mean value of $B_{k,1}(x,y)$ on the unit sphere. The parallel product formula for $a=2$ (the Dunkl case) was given by R\"osler in \cite{R2}.
\end{rem}

In the proof of the positivity of the spherical mean operators for $a=1$ in Theorem \ref{positivity}, the classical method for $a=2$ in \cite{R2} to prove firstly for Schwartz functions fails  because the Schwartz space $\mathcal S(\mathbb R^N)$ is not invariant under the $(k,1)$-generalized Fourier transform $F_{k,1}$ (see \cite{GIT3}), contrary to the claim in \cite{Jo}. 
Before the claim in \cite{Jo} was known to be false, a proof of Theorem \ref{positivity} (1) was given in \cite{Iv1}. But we must prove without reliance on the false conclusion.
We will first prove the positivity for functions in the space $W_{k,1}(\mathbb R^N)$, which was introduced in \cite[Section 3]{BSKO} and is invariant under the generalized Fourier transform $F_{k,1}$.
And then we can extend the positivity to all functions $f\in C_0(\mathbb{R}^N\backslash{\left\{0\right\}})$ by a density argument because it will be shown in Section 2 that the space $W_{k,1}(\mathbb R^N)$ invariant under $F_{k,1}$ is ${\left\|\cdot\right\|}_\infty$-dense in $C_0(\mathbb{R}^N\backslash{\left\{0\right\}})$.

The function $d\left(x,y\right)=\sqrt{\left|x\right|+\left|y\right|-\sqrt{2\left(\left|x\right|\left|y\right|+\left\langle x,y\right\rangle\right)}}$ (see also Appendix), which was shown in \cite{Te3} to be the metric corresponding to the setting of $(k,1)$-generalized Fourier analysis, was derived from the structure of the radial formula of the $(k,1)$-generalized translation in \cite{BD}
\begin{align}\label{radial}\tau_yf(x)=&\frac{\Gamma\left(\frac{N-1}2+\left\langle k\right\rangle\right)}{\sqrt\pi\Gamma\left(\frac{N-2}2+\left\langle k\right\rangle\right)}\times\notag\\&V_k\left(\int_{-1}^1f_0\left(\left|x\right|+\left|y\right|-\sqrt{2\left(\left|x\right|\left|y\right|+\left\langle\cdot,y\right\rangle\right)}u\right)\left(1-u^2\right)^{\frac N2+\left\langle k\right\rangle-2}du\right)\left(x\right),\end{align}
where $2\langle
k\rangle+N-2>0$, $f(x)=f_0\left(\left|x\right|\right)\in \mathcal L_k^1\left(\mathbb{R}^N\right)$ and
$\mathcal L_k^1\left(\mathbb{R}^N\right):=\{f\in
L^1\left(\mathbb{R}^N,\vartheta_{k,1}\left(x\right)dx\right):F_{k,1}\left(f\right)\in
L^1\left(\mathbb{R}^N,\vartheta_{k,1}\left(x\right)dx\right)\}.$ 
Note that \eqref{radial} holds true for radial Schwartz functions since $F_{k,1}(\mathcal S_{rad}(\mathbb{R}^N))$ consists of rapidly decreasing functions
at infinity, according to the Proposition 5.5.(ii) in \cite{GIT3} for even functions on the real line.
And \eqref{radial} can be extended to all continuous radial
functions on ${{L}^{2}}\left(
{{\mathbb{R}}^{N}},{\vartheta_{k,1}}\left( x \right)dx \right)$
from a similar argument to the Lemma 3.4 in \cite{DW}.

The Schwartz space $\mathcal S(\mathbb{R}^N)$ is not invariant under the $(k,1)$-generalized Fourier transform $F_{k,1}$ because the $(k,1)$-generalized Fourier transform of a nontrivial function may not be differentiable at the origin. 
And when $N\geq 2$, it is unknown whether all the functions in $F_{k,1}\left(\mathcal S(\mathbb{R}^N)\right)$ are rapidly decreasing at infinity (see \cite {GIT3} for the one dimensional case).
We shall consider the following invariant subspace of the $(k,1)$-generalized Fourier transform
$$\widehat{\mathcal S}(\mathbb{R}^N):=\left\{f\;\mathrm{decays}\;\mathrm{faster}\;\mathrm{than}\;\mathrm{any}\;\mathrm{polynomial}:\;F_{k,1}(f)\;\mathrm{decays}\;\mathrm{faster}\;\mathrm{than}\;\mathrm{any}\;\mathrm{polynomial}\right\}.$$
Here we say that $f$ decays faster than any polynomial if for all $m\in\mathbb{N}$, there exists $C_m>0$, such that $\left|f(x)\right|\leq\frac{C_m}{{(1+\left|x\right|)}^m}$.
Obviously, $W_{k,1}(\mathbb{R}^N)\subset \widehat{\mathcal S}(\mathbb{R}^N)\subset C^\infty(\mathbb{R}^N\backslash\{0\})$. And it can be observed that $M_f\in C^\infty(\mathbb{R}^N\backslash\{0\}\times \mathbb R_+)$ if $f\in
\widehat{\mathcal S}(\mathbb{R}^N)$, which allows us to consider related problems in partial differential equations.
There are also several attempts to obtain an invariant subspace  under the $(k,a)$-generalized Fourier transform (see e.g. \cite{Iv} and \cite{FN}).
The following theorem is parallel to the Proposition 5.3 in \cite{Me1}.
\begin{thm}\label{wave equation}
Assume $2\langle k\rangle+N-2>0$. 
Let $f\in
\widehat{\mathcal S}(\mathbb{R}^N)$. Then
$\mathcal R_{\lambda_1}^{-1}M_f(\cdot,\frac12t^2)$, where $\mathcal R_\alpha$ is the Riemann-Liouville transform and $\lambda_1=2\left\langle k\right\rangle+N-2$, is the unique solution to the following deformed wave equation in $C(\mathbb{R}^N)\cap C^\infty(\mathbb{R}^N\backslash\{0\})\cap L^2\left(\mathbb{R}^N,\vartheta_{k,1}\left(x\right)dx\right)$
\begin{align}\label{wave}
 &u_{tt}-2\left|x\right|\triangle_ku=0\;\;(x\in\mathbb{R}^N\backslash\{0\},\;t>0);\notag\\
 &u(x,0) = f(x), \,\, u_t (x,0) = 0, \;\;x\in \mathbb R^N.
\end{align}
\end{thm}

We can then deduce a Huygens' principle from the support of the representing measure $\sigma_{x,t}^{k,1}$ of $M_f(x,t)$ in Theorem \ref{positivity} (2).
\begin{thm}\label{domain of dependence}
(Weak Huygens' principle) Assume $2\langle k\rangle+N-2>0$. Then the solution to the Cauchy problem \eqref{wave} at a given point $(x,t)\in\mathbb{R}^N\backslash\{0\}\times
\mathbb R_+$ depends only on the values of the initial data $f$ in the union 
$$\bigcup_{g\in G}\{\xi\in \mathbb R^N: \sqrt 2d(\xi,gx)\leq  t\}.$$ 
\end{thm}

This paper is organized as follows. In Section 2, we review some results in Dunkl theory (ref. \cite{Du, Du1,Du2}) and the $(k,a)$-generalized Fourier analysis developed in \cite{BSKO} as necessary tools for proving the main results. In Section 3, we will show that $W_{k,1}\left(\mathbb{R}^N\right)$ is ${\left\|\cdot\right\|}_\infty$-dense in  $C_0(\mathbb{R}^N\backslash{\left\{0\right\}})$, as a preparatory work before showing the positivity of the generalized spherical mean operator. In Section 4 we will prove the first part of Theorem \ref{positivity}. 
The equivalence to the positive radial product formula \eqref{product formula2} will be derived in
Section 5. And Section 5 will be devoted to the analysis of the representing measure $\sigma_{x,t}^{k,1}$ to complete the proof of the second and third part of Theorem \ref{positivity}. We will then study the deformed wave equation \eqref{wave} and prove Theorem \ref{wave equation} and Theorem \ref{domain of dependence} in Section 5.
In Section 6 we consider the solution to the deformed wave equation and write the solution in terms of the generalized spherical mean operator.
In the final appendix we show that the metric $\sqrt 2d(x,y)$  for $N\geq 2$ is in fact a Riemannian distance, as a supplementary, linking $(k,1)$-generalized Fourier analysis with Riemannian geometry.
We denote $\lambda_{k,a,m}:=\frac{2m+2\left\langle
k\right\rangle+N-2}a$ and $\lambda_a:=\frac{2\left\langle k\right\rangle+N-2}a$ in this paper. The function classes $C_c^\infty$, $C_0^\infty$ and $C^\infty$ are defined in the standard manner.

\section{Preliminaries}
\subsection{Dunkl theory}

\

Dunkl theory is a far-reaching generalization of Fourier analysis and special function theory about root system with a rich structure parallel to ordinary Fourier analysis. 
In the remainder of this subsection we will review the framework of Dunkl theory to introduce the essential tools we will need.

Given a (reduced but not necessarily  crystallographic) root system $R$ in the Euclidean space $\mathbb R^N$, denote by $G$ the
finite subgroup of the orthogonal group $O(N)$ generated by the reflections
$\sigma_\alpha$ associated to the root system. Define a \textit{multiplicity function}
$k:R\rightarrow \mathbb C$ such that $k$ is $G$-invariant, that is,
$k\left(\alpha\right)=k\left(\beta\right)$ if $\sigma_\alpha$ and
$\sigma_\beta$ are conjugate. We denote $R^+$ to be any fixed positive subsystem of $R$ and $\left\langle
k\right\rangle:=\sum_{\alpha\in R^+}k(\alpha)$. 
Assume the root system $R$ is normalized (i.e., $\left\langle\alpha,\alpha\right\rangle=2$ for all vectors $\alpha\in R$) without loss of generality.
In \cite{Du}, C. F. Dunkl constructed a kind of
 differential difference operator as follows:
\begin{align}\label{reflectioninvariant}\triangle_k f(x)=\triangle f\left(x\right)+2\sum_{\alpha\in
R^+}k\left(\alpha\right)\frac{\left\langle\nabla
f,\alpha\right\rangle}{\left\langle\alpha,x\right\rangle}-2\sum_{\alpha\in
R^+}k\left(\alpha\right)\frac{f\left(x\right)-f\left(\sigma_\alpha\left(x\right)\right)}{\left\langle\alpha,x\right\rangle^2},\end{align}
where
$\nabla$ is the Euclidean gradient and $\triangle$ is the Euclidean Laplacian. The operators commute with the action of the finite reflection group $G$ and are symmetric on the Hilbert space $L^2\left(\mathbb{R}^N,\vartheta_k\left(x\right)dx\right)$, where 
$$\vartheta_k\left(x\right)=\prod_{\alpha\in R^+}\vert\left\langle\alpha,\;x\right\rangle\vert^{2k(\alpha)}.$$
Let $P_m$ be the space of homogeneous
polynomials on $\mathbb{R}^N$ of degree $m$.
It was shown by Dunkl \cite {Du} that the restrictions $\mathcal
H_k^m\left(\mathbb{R}^N\right)\vert_{\mathbb S^{N-1}},\;m=0,1,2,\dots,$ of the spaces $\mathcal H_k^m\left(\mathbb{R}^N\right):=P_m\cap ker\triangle_k,\;m=0,1,\cdots,$ to the unit sphere $\mathbb S^{N-1}$ are orthogonal to each other with respect to  $\vartheta_k\left(x\right)d\sigma$, where $d\sigma$ is the spherical measure. The spaces $\mathcal
H_k^m\left(\mathbb{R}^N\right)\vert_{\mathbb S^{N-1}}$ are called spherical $k$-harmonics and the operator $\triangle_k$ is called \textit{Dunkl Laplacian}. The following is the
spherical harmonics decomposition in Dunkl setting
\begin{equation}\label{decomposition}L^2\left(\mathbb S^{N-1},\vartheta_k\left(x'\right)d\sigma(x')\right)=\sum_{m\in \mathbb N}^\oplus\mathcal H_k^m\left(\mathbb{R}^N\right)\vert_{\mathbb S^{N-1}},\;x'\in \mathbb S^{N-1}.\end{equation}
The eigenfunction $E_k(\cdot,-iy)$ of the Dunkl Laplacian $\triangle_k$ for fixed $y$ is the integral kernel of the generalized Fourier transform $F_k$ called Dunkl transform. 

The classical Lapalcian $\triangle$ and the Dunkl Laplacian $\triangle_k$ are intertwined by a
Laplace-type operator (see \cite{Du2}) which is a homeomorphism on $C^\infty(\mathbb{R}^N)$
\begin{align}V_kf(x)=\int_{\mathbb{R}^N}f(y)d\mu_x(y),\end{align}
that is, 
$\triangle_k\circ V_k=V_k\circ\triangle.$
It is associated to a family of probability measures
$\left\{\mu_x\vert\;x\in\mathbb{R}^N\right\}$ with compact support
(see \cite{R2}).
Specifically, the support of $\mu_x$ is contained in the convex
hull $co(G.x)$, where $G.x=\left\{g.x\vert\;g\in
G\right\}$ is the orbit of $x$. For any Borel set $B$ and any
$r>0$, $g\in G$, the probability measures satisfy
$$\mu_{rx}\left(B\right)=\mu_x\left(r^{-1}B\right),\;\mu_{gx}\left(B\right)=\mu_x\left(g^{-1}B\right).$$

\subsection{ The $(k,a)$-generalized Fourier transform}

\

In the following subsection we introduce some definitions and results in the development of the $(k,a)$-generalized Fourier analysis initiated in S. Ben Sa\"id, T. Kobayashi and B. \O rsted \cite{BSKO}.

Consider the weight function ${\vartheta_{k,a}}\left( x
\right)={{\left| x \right|}^{a-2}}\vartheta_{k}(x)$. It reduces to
$\vartheta_k(x)$ when $a=2$ and also on the unit sphere $\mathbb S^{N-1}$. 
Let
$${{\Delta }_{k,a}}={{\left| x \right|}^{2-a}}{{\Delta }_{k}}-{{\left| x \right|}^{a}},a>0.$$
It is a symmetric operator on ${{L}^{2}}\left(
{{\mathbb{R}}^{N}},{\vartheta_{k,a}}\left( x \right)dx \right)$ with only negative discrete spectrum.

For the polar coordinates $x=rx'(r>0,\;x'\in
\mathbb S^{N-1})$,
we have the following
unitary isomorphism from the spherical harmonic decomposition \eqref{decomposition} of
$L^2\left(\mathbb S^{N-1},\vartheta_k\left(x'\right)d\sigma(x')\right)$ (see \cite[(3.25)]{BSKO})
$$\sum_{m\in \mathbb N}^\oplus(\mathcal H_k^m\left(\mathbb{R}^N\right){\vert_{\mathbb S^{N-1}})\otimes L^2}\left({\mathbb{R}}_+,r^{\;2\left\langle k\right\rangle+N+a-3}dr\right)\xrightarrow\sim L^2\left(\mathbb{R}^N,\vartheta_{k,a}\left(x\right)dx\right).$$
Define the Laguerre polynomial as
$$L_l^\mu(t):=\sum_{j=0}^l\frac{{(-1)}^j\Gamma(\mu+l+1)}{(l-j)!\Gamma(\mu+j+1)}\frac{t^j}{j!},\;\mathrm{Re}\mu>-1.$$
And
consider the following linear operator 
$S_a: C^{\infty}\left(\mathbb S^{N-1}\right) \otimes C^{\infty}\left(\mathbb{R}_{+}\right) \rightarrow C^{\infty}\left(\mathbb{R}^N \backslash\{0\}\right)$, where $S_a$ is given as
$$
S_a(p\otimes g)(x):=p(x')\exp\left(-\frac1ar^a\right)g\left(\frac2ar^a\right).
$$
For $l,m\in\mathbb{N}$, and $p\in\mathcal
H_k^m\left(\mathbb{R}^N\right)\vert_{\mathbb S^{N-1}}$, we introduce the following functions on $\mathbb R^N$
\begin{align*}\mathrm\Phi_l^{(a)}(p,x):&=S_a\left(p\otimes L_l^{\left(\lambda_{k,a,m}\right)}\right)(x)\\&=p(x')r^mL_l^{\lambda_{k,a,m}}\left(\frac2ar^a\right)\exp\left(-\frac1ar^a\right).\end{align*}
Denote 
$$W_{k,a}\left(\mathbb{R}^N\right):=\mathbb{C}\text{-}span\left\{\left.\mathrm\Phi_{\mathit l}^{\mathit(\mathit a\mathit)}\mathit{(p,\cdot)}\right|l\in\mathbb{N},\;m\in\mathbb{N},\;p\in\mathcal H_k^m\left(\mathbb{R}^N\right)\right\}.$$
It was shown in \cite{BSKO} that $W_{k,a}\left(\mathbb{R}^N\right)$ is a dense subspace of $L^2\left(\mathbb{R}^N,\vartheta_{k,a}\left(x\right)dx\right)$.
And the operator ${{\Delta }_{k,a}}$ was defined on $L^2\left(\mathbb{R}^N,\vartheta_{k,a}\left(x\right)dx\right)$ as an unbounded operator with the dense domain $W_{k,a}\left(\mathbb{R}^N\right)$.

The $(k,a)$-generalized Laguerre semigroup ${\mathcal I}_{k,a}\left(z\right)$ was then defined as follows in \cite{BSKO} with the infinitesimal generator $\frac1a\triangle_{k,a}$, 
$${\mathcal I}_{k,a}\left(z\right):=\exp \left(\frac za\triangle_{k,a}\right),\;\Re z\geq0.$$
When taking the boundary value $z=\frac{\pi i}2$, the semigroup  $\mathcal I_{k,a}\left(z\right)$ reduces to the $(k,a)$-generalized Fourier transform $F_{k,a}$, i.e.,\\
\begin{align}\label{Fka}F_{k,a}=c{\mathcal I}_{k,a}\left(\frac{\pi i}2\right),
\end{align}
where $c=e^{i\pi(\frac{2\left\langle k\right\rangle+N+a-2}{2a})}.$
The generalized Fourier transform includes
the Fourier transform ($k\equiv0$ and $a=2$), the Kobayashi-Mano Hankel transform \cite{KM1, KM2} ($k\equiv0$ and $a=1$), and the Dunkl transform \cite{Du2} ($k\geq0$ and $a=2$).
For $l, m \in \mathbb{N}$ and $p \in \mathcal{H}_k^m\left(\mathbb{R}^N\right)$, $\Phi_{l}^{(a)}(p, \cdot)$ is an eigenfunction of $F_{k, a}$, i.e.,
$$
F_{k, a}\left(\Phi_{l}^{(a)}(p, \cdot)\right)=e^{-i \pi\left(l+\frac{m}{a}\right)} \Phi_{l}^{(a)}(p, \cdot).
$$
Therefore, $W_{k,a}\left(\mathbb{R}^N\right)$ is an invariant subspace of the $(k,a)$-generalized Fourier transform $F_{k,a
}$. And we have the Plancherel formula of the $(k,a)$-generalized Fourier transform, i.e., ${\left\|F_{k,a}f\right\|}_{L^2}={\left\|f\right\|}_{L^2}$.

By Schwartz kernel theorem, the $(k,a)$-generalized Fourier
transform has the following integral representation on ${{L}^{2}}\left(
{{\mathbb{R}}^{N}},{\vartheta_{k,a}}\left( x \right)dx \right)$ (see \cite[(5.8)]{BSKO}) 
\begin{align}\label{Fourier transform}F_{k,a}f\left(\xi\right)=c_{k,a}\int_{\mathbb{R}^N}f\left(y\right)B_{k,a}\left(\xi,y\right)\vartheta_{k,a}\left(y\right)dy,\;\;\xi\in\mathbb{R}^N,\end{align}
where $c_{k,a}=\left(\int_{\mathbb{R}^N}\exp\left(-\left|x\right|^a\right)\vartheta_{k,a}\left(x\right)dx\right)^{-1}$ and $B_{k,a}\left(x,y\right)$ is a symmetric kernel. Moreover,
\begin{align}\label{kernelproperty}
B_{k,a}\left(\alpha x,y\right)=B_{k,a}\left(x,\alpha y\right)\;\mathrm{and}\;B_{k,a}\left(gx,gy\right)=B_{k,a}\left(x,y\right)
\end{align}
for all $x,y\in\mathbb{R}^N$, $\alpha\in\mathbb{R}$ and $g\in G$.
The integral kernel $B_{k,a}(x,y)$ of the $(k,a)$-generalized Fourier transform takes the place of the exponential function $e^{-i\left\langle x,y\right\rangle}$ in classical Fourier transform. It is the eigenfunction of the operator $\left|x\right|^{2-a}\triangle_k$ for any fixed $y$ (see \cite[Theorem 5.7]{BSKO}), i.e.,
$$\left|x\right|^{2-a} \Delta_k^x B_{k,a}(x,y)= -\left|y\right|^a B_{k,a}(x,y).$$
So, we can consider the operator $\left|x\right|^{2-a} \Delta_k$ as the $a$-deformed Dunkl Laplacian in $(k,a)$-generalized Fourier analysis.

For $a={\textstyle\frac2{n}},\;n\in{\mathbb{N}}_+$, $2\left\langle k\right\rangle+N+a-2>0$,
we have the inversion formulas of the $(k,a)$-generalized
Fourier transform (see \cite[Theorem 5.3]{BSKO}), i.e.,
\begin{align}\label{inversion}&\left(F_{k,a}\right)^{-1}=F_{k,a},\;\;\;\;\;\;\;\;\;\;\;\;\;\;\;\;\;\;\;\;\;\mathrm{if}\;a={\textstyle\frac1r},\;r\in{\mathbb{N}}_+,\\&\left(F_{k,a}^{-1}f\right)\left(x\right)=\left(F_{k,a}f\right)\left(-x\right),\;\;\mathrm{if}\;a={\textstyle\frac2{2r+1}},\;r\in\mathbb{N}.\notag\end{align}
And for $N\geq2,\;a={\textstyle\frac2{n}},\;n\in{\mathbb{N}}_+$, the integral kernel $B_{k,a}(x,y)$ is uniformly bounded by $1$ (see \cite{CD}), i.e.,
\begin{align*}\left|B_{k,a}\left(x,y\right)\right|\leq\left|B_{k,a}\left(0,y\right)\right|=1.\end{align*} 
It still remains an open problem on more general conditions under which the integral kernel $B_{k,a}(x,y)$ is uniformly bounded.
For the one dimensional case the the necessary and sufficient condition is
$4k+a-2\geq0$, $k\geq 0$ and $a>0$ (see \cite[Proposition 2.1]{GIT3}).
In the cases that $B_{k,a}(x,y)$ is uniformly bounded, the integral representation \eqref{Fourier transform} of the $(k,a)$-generalized Fourier transform can be extended to all functions in ${{L}^{1}}\left(
{{\mathbb{R}}^{N}},{\vartheta_{k,a}}\left( x \right)dx \right)$.

Assume the conditions such that $B_{k,a}(x,y)$ is uniformly bounded. For $a={\textstyle\frac2{n}},\;n\in{\mathbb{N}}_+$, $2\left\langle k\right\rangle+N+a-2>0$, one can define the
$(k,a)$-generalized translations $\tau_y$ on ${{L}^{2}}\left(
{{\mathbb{R}}^{N}},{\vartheta_{k,a}}\left( x \right)dx \right)$ as
\begin{align}\label{deftranslation}
F_{k,a}\left(\tau_yf\right)\left(\xi\right)=B_{k,a}\left({(-1)}^ny,\xi\right)F_{k,a}\left(f\right)\left(\xi\right),\;\;\xi\in\mathbb{R}^N.
\end{align}
The above definition makes sense because $F_{k,a}$ is an isometry on
${{L}^{2}}( {{\mathbb{R}}^{N}},$ ${\vartheta_{k,a}}\left( x
\right)dx )$.
In this case for $f\in \mathcal L_k^1\left(\mathbb{R}^N\right)$, where
$\mathcal L_k^1\left(\mathbb{R}^N\right):=\{f\in
L^1\left(\mathbb{R}^N,\vartheta_{k,a}\left(x\right)dx\right):F_{k,a}\left(f\right)\in
L^1\left(\mathbb{R}^N,\vartheta_{k,a}\left(x\right)dx\right)\},$ 
the $(k,a)$-generalized translations
can also be written via integrals as
\begin{align}\label{translation2} \tau_yf(x)=c_{k,a}\int_{\mathbb{R}^N}B_{k,a}\left((-1)^nx,\xi\right)B_{k,a}\left((-1)^ny,\xi\right)F_{k,a}\left(f\right)\left(\xi\right)\vartheta_{k,a}\left(\xi\right)d\xi,\;\;\mathrm{if}\;a={\textstyle\frac2{n}},\;n\in{\mathbb{N}}.\end{align}
Since $\widehat{\mathcal S}(\mathbb{R}^N)$ is a subspace of $\mathcal L_k^1\left(\mathbb{R}^N\right)$, the formulas \eqref{translation2} hold true for all $f\in \widehat{\mathcal S}(\mathbb{R}^N)$.
The $(k,a)$-generalized translations satisfy the following properties:\\
(1). For every $x,y\in \mathbb R^N$,
\begin{align}\label{prop1}\tau_yf(x)=\tau_xf(y),\;f\in\mathcal L_k^1\left(\mathbb{R}^N\right).\end{align} (2). For every $y\in
\mathbb R^N$, and $f\in\mathcal L_k^1\left(\mathbb{R}^N\right)$, $g\in L^1\left(\mathbb{R}^N,\vartheta_{k,a}\left(x\right)dx\right)\cap L^\infty\left(\mathbb{R}^N,\vartheta_{k,a}\left(x\right)dx\right)$, 
\begin{align*}\int_{\mathbb{R}^N}\tau_yf\left(x\right)g\left(x\right)\vartheta_{k,a}\left(x\right)dx=\int_{\mathbb{R}^N}f\left(x\right)\tau_yg\left(x\right)\vartheta_{k,a}\left(x\right)dx,\;\;\;\;\;\;\;\mathrm{if}\;a={\textstyle\frac1r},\;r\in{\mathbb{N}}_+, \end{align*}

\;\;\;and 
\begin{align*}\int_{\mathbb{R}^N}\tau_yf\left(x\right)g\left(x\right)\vartheta_{k,a}\left(x\right)dx=\int_{\mathbb{R}^N}f\left(x\right)\tau_{-y}g\left(x\right)\vartheta_{k,a}\left(x\right)dx,\;\;\;\;\;\;\;\mathrm{if}\;a={\textstyle\frac2{2r+1}},\;r\in{\mathbb{N}}. \end{align*}
The Property (2) allows us to define the $(k,a)$-generalized translations on $L^p\left(\mathbb{R}^N,\vartheta_{k,a}\left(x\right)dx\right),\;1\leq p\leq\infty$ in the distributional sense.
In the following we denote $f(x\ast y):=\tau_xf(y)$ for convenience in view of the Property (1).
We can extend such definition and define
$$B_{k,a}(x\ast y,z):=B_{k,a}(x,z)B_{k,a}(y,z)$$
from the symmetry property of the $(k,a)$-generalized Fourier transform $F_{k,a}$.

In the cases that $B_{k,a}\left(x,y\right)$ is uniformly bounded,
we shall consider also the $(k.a)$-generalized Fourier transform on a bounded Borel measure $\mu$ on $\mathbb R^N$
$$F_{k,a}\left[\mu\right](\xi):=\int_{\mathbb{R}^N}B_{k,a}(\xi,y)d\mu(y),\quad\xi\in\mathbb{R}^N.$$
Then from Fubini theorem we notice
\begin{align}\label{Borel measure}
\int_{\mathbb{R}^N}f(x)F_{k,a}\left[\mu\right](x)\vartheta_{k,a}(x)dx=\int_{\mathbb{R}^N}F_{k,a}(f)d\mu.
\end{align}
for $f\in {L^1}\left(
{{\mathbb{R}}^{N}},{\vartheta_{k,a}}\left( x \right)dx \right)$.
And according to \eqref{Borel measure} the $(k,a)$-generalized Fourier transform is injective on bounded Borel measures, i.e., if $F_{k,a}\left[\mu\right]=0$, then $\mu=0$, since $F_{k,a}\left(L^1\left(\mathbb{R}^N,\vartheta_{k,a}\left(x\right)dx\right)\right)$ is a  ${\left\|\cdot\right\|}_\infty$-dense subspace in $C_0(\mathbb{R}^N)$ by the locally compact version of Stone--Weierstrass theorem.

\section{${\left\|\cdot\right\|}_\infty$-Density of the space $W_{k,1}\left(\mathbb{R}^N\right)$ in $C_0(\mathbb{R}^N\backslash{\left\{0\right\}})$}

Before we can prove the main results, we need to show that the space $W_{k,1}\left(\mathbb{R}^N\right)$ is dense in $C_0(\mathbb{R}^N\backslash{\left\{0\right\}})$ (every function $f\in C_0(\mathbb{R}^N\backslash{\left\{0\right\}})$ can be approximated by the functions in $W_{k,1}\left(\mathbb{R}^N\right)$ with respect to ${\left\|\cdot\right\|}_\infty$), so that we need only to prove the positivity-preserving property on $W_{k,1}\left(\mathbb{R}^N\right)$.
In the following we denote $\psi_{l,m}(r):=\psi_{l,m}^{(1)}(r)$ and $\lambda_{k,m}:=\lambda_{k,1,m}$.

For fixed $m\in\mathbb{N}$ and a multiplicity function $k$
satisfying $\lambda_{k,m}>-1$ we set
\begin{equation}\label{psi}\psi_{l,m}(r):=\left(\frac{2^{\lambda_{k,m}+1}\Gamma(l+1)}{\Gamma(\lambda_{k,m}+l+1)}\right)^{1/2}r^mL_l^{\lambda_{k,m}}\left(2r\right)\exp\left(-r\right).\end{equation}
It was shown in \cite[Proposition 3.15]{BSKO} that $\left\{\psi_{l,m}(r):\;l\in\mathbb{N}\right\}$ forms
an orthonormal basis in
$L^2({\mathbb{R}}_+,r^{\;2\left\langle
k\right\rangle+N-2}$ $dr)$.
We will show that the functions in $C_0^\infty(\mathbb{R}_+)$ can also be approximated by the linear combinations of $\left\{\psi_{l,m}(r):\;l\in\mathbb{N}\right\}$.
For $m=0$ this was already shown in \cite{No}.

\begin{prop}\label{psilm}
For a fixed $m\in\mathbb{N}$, $\mathbb C$-$span\left\{\psi_{l,m}(r):\;l\in\mathbb{N}\right\}$ is dense in $C_0^\infty(\mathbb{R}_+)$ with respect to the norm ${\left\|\cdot\right\|}_\infty$.
\end{prop}
\begin{proof}
It suffices to prove that every function $f\in C_c^\infty\left(\mathbb{R}_+\right)$ can be approximated by the linear combinations of $\left\{\psi_{l,m}(r):\;l\in\mathbb{N}\right\}$.
We first need the following estimate of $\psi_{l,m}(r)$
\begin{align}\label{estimate}\left|\psi_{l,m}(r)\right|\leq\left\{\begin{array}{lc}c{(2l+\lambda_{k,m}+1)}^{\lambda_{k,m}},&0<r\leq{\textstyle\frac32}(2l+\lambda_{k,m}+1),\\c\exp(-\gamma r),&r>{\textstyle\frac32}{\textstyle(}{\textstyle2}{\textstyle l}{\textstyle+}{\textstyle{\lambda}_{k,m}}{\textstyle+}{\textstyle1}{\textstyle)}{\textstyle,}\end{array}\right.\end{align}
where $c$ and $\gamma$ are constants independent of $k$ and $r$. 
This estimate (see c.f. \cite{No}) is a consequence of Muckenhoupt's \cite{Mu2} estimate of Laguerre functions.
For a fixed function $f\in C_c^\infty\left(\mathbb{R}_+\right)$, define
$$S_Lf=\sum_{l=0}^L{\left\langle\psi_{l,m},f\right\rangle}_{L^2\left({\mathbb{R}}_+,r^{\;2\left\langle k\right\rangle+N-2}dr\right)}\psi_{l,m}.$$ 
Then $S_Lf\rightarrow f$ in $L^2\left({\mathbb{R}}_+,r^{\;2\left\langle k\right\rangle+N-2}dr\right)$. And there exists a subsequence $S_{L_k}f$ such that $S_{L_k}f\rightarrow f,\;a.e.$
It remains to show that $S_{L}f$ is uniformly fundamental.
In \cite{Te2} the author defined the deformed Laguerre operator
\begin{equation}\label{Laguerre}
L_{1,\alpha}=-r\frac{d^2}{dr^2}+r-\left(\alpha+1\right)\frac d{dr}
\end{equation}
and showed that the Laguerre functions $\widetilde{\varphi}_l^{1,\lambda_{k,m}}(r):=r^{-m}\psi_{l,m}(r)$ are eigenfunctions of the Laguerre operator $L_{1,\lambda_{k,m}}$, i.e., 
$$L_{1,\lambda_{k,m}}\widetilde{\varphi}_l^{1,\lambda_{k,m}}=\left(2l+\lambda_{k,m}+1\right)\widetilde{\varphi}_l^{1,\lambda_{k,m}},\;\;l=0,1,\cdots.$$
For the representation  $\Omega_{k,a}^{(m)}$ of $S\widetilde{L(2,}\mathbb{R})$ on $L^2\left({\mathbb{R}}_+,r^{\;2\left\langle k\right\rangle+N-2}dr\right)$ (see \cite[Section 4]{BSKO} for the detailed definition of $\Omega_{k,a}^{(m)}\left(\mathrm{Exp}\left(-\mathrm z\mathbf k\right)\right)$, $\Re\;z\geq0$), the derivative of $\Omega_{k,a}^{(m)}\left(\mathrm{Exp}\left(-\mathrm z\mathbf k\right)\right)$ can be expressed as the following symmetric operator via the deformed Laguerre operator (see \cite [Theorem 5.1]{Te2})
$$\mathrm d\Omega_{k,1}^{(m)}(-\mathbf k)f(r)=-r^mL_{1,\lambda_{k,m}}\left(\left(\cdot\right)^{-m}f\right)\left(r\right).$$
Therefore,
for any $n\in\mathbb{N}$,
\begin{align}\label{symmetry}{\left\langle\psi_{l,m},f\right\rangle}_{L^2\left({\mathbb{R}}_+,r^{\;2\left\langle k\right\rangle+N-2}dr\right)}&\notag={\left\langle\left(\mathrm d\Omega_{k,a}^{(m)}(\mathbf k)\right)^{-n}\psi_{l,m},\left(\mathrm d\Omega_{k,a}^{(m)}(\mathbf k)\right)^nf\right\rangle}_{L^2\left({\mathbb{R}}_+,r^{\;2\left\langle k\right\rangle+N-2}dr\right)}\\&=\left(2l+\lambda_{k,m}+1\right)^{-n}{\left\langle\psi_{l,m}(r),\left(\mathrm d\Omega_{k,a}^{(m)}(\mathbf k)\right)^nf\right\rangle}_{L^2\left({\mathbb{R}}_+,r^{\;2\left\langle k\right\rangle+N-2}dr\right)}.\end{align}
Let $1\leq N_1<N_2$. Then from \eqref{symmetry}, \eqref{estimate} and Cauchy--Schwartz inequality we have
\begin{align*}\left|S_{N_1}f-S_{N_2}f\right|&\leq\sum_{l=N_1}^{N_2}\left|\left\langle\psi_{l,m},f\right\rangle\right|\left|\psi_{l,m}\right|\\&\leq\sum_{l=N_1}^{N_2}c\left(2l+\lambda_{k,m}+1\right)^{-n}\left|\left\langle\psi_{l,m}(r),\left(\mathrm d\Omega_{k,a}^{(m)}(\mathbf k)\right)^nf(r)\right\rangle\right|\left(2l+\lambda_{k,m}+1\right)^{\lambda_{k,m}}\\&\leq\sum_{l=N_1}^{N_2}c\left(2l+\lambda_{k,m}+1\right)^{-n}{\left\|\left(\mathrm d\Omega_{k,a}^{(m)}(\mathbf k)\right)^nf\right\|}_{L^2}\left(2l+\lambda_{k,m}+1\right)^{\lambda_{k,m}}.\end{align*}
If we take $n>\lambda_{k,m}$, then the last expression tends to $0$ as $N_1,\;N_2\rightarrow\infty$.
\end{proof}

\begin{thm}
$W_{k,1}\left(\mathbb{R}^N\right)=\mathbb{C}\text{-}span\left\{\left.\mathrm\Phi_{\mathit l}^{( 1)}\mathit{(p,\cdot)}\right|l\in\mathbb{N},\;m\in\mathbb{N},\;p\in\mathcal H_k^m\left(\mathbb{R}^N\right)\right\}$ is dense in $C_0(\mathbb{R}^N\backslash{\left\{0\right\}})$ with respect to the norm ${\left\|\cdot\right\|}_\infty$.
\end{thm}

\begin{proof}
From the Weierstrass approximation theorem,
the linear span of the polynomials in the union of the spaces $\mathcal
H_k^m\left(\mathbb{R}^N\right)\vert_{\mathbb S^{N-1}},\;m=0,1,2,\dots$ is a dense subspace of $C^\infty(\mathbb S^{N-1})$.
We can then get that every function in 
$C_0^\infty(\mathbb{R}^N\backslash{\left\{0\right\}})$
can be approximated by the elements in $W_{k,1}\left(\mathbb{R}^N\right)$ with respect to the norm ${\left\|\cdot\right\|}_\infty$ combining Proposition \ref{psilm}.
And therefore $W_{k,1}\left(\mathbb{R}^N\right)$ is ${\left\|\cdot\right\|}_\infty$-dense in $C_0(\mathbb{R}^N\backslash{\left\{0\right\}})$, since $C_0^\infty(\mathbb{R}^N\backslash{\left\{0\right\}})$ is a ${\left\|\cdot\right\|}_\infty$-dense subspace of $C_0(\mathbb{R}^N\backslash{\left\{0\right\}})$.
\end{proof}

\section{The positivity of the spherical mean operator}

Now we are able to investigate the $(k,a)$-generalized spherical mean operator and show that the operator is positivity-preserving for $a=1$ after the preparatory work in the last section. 
We define the generalized spherical mean operator $f\mapsto M_f$ on $C_0(\mathbb{R}^N)$ as
\begin{align}\label{spherical mean} M_f(x,t):= \frac{1}{d_{k,a}} \int_{\mathbb S^{N-1}} f(x*ty) {\vartheta_{k,a}}(y) d\sigma(y),
 \quad (x\in \mathbb R^N, \, t\geq 0),\end{align}
where $d\sigma$ is the spherical measure and $d_{k,a} = \int_{\mathbb S^{N-1}} {\vartheta_{k,a}}\left( x
\right)d\sigma(x)$. For $a=2$, it was shown by R\"osler \cite{R2} that the operator is positivity-preserving  and is uniquely represented by a continuous and compactly supported probability measure $\sigma_{x,t}^{k,2}(\xi)$.
We will show the positivity-preserving property of the operator for $a=1$ in this section.

Firstly, we consider the one dimensional case assuming $ k\geq1/2$. 
For the one dimensional case it was already shown in \cite{Iv} for $a={\textstyle\frac1r},\;r\in{\mathbb{N}}_+$, $ r(2k-1)>-1/2$.
We will give the proof for all $a={\textstyle\frac2n},\;n\in{\mathbb{N}}_+$, $k\geq1/2$ here. For $N=1$, $a={\textstyle\frac2n},\;n\in{\mathbb{N}}_+$, and $f=f_e+f_o$ written as a sum of even and odd functions,  H. Mejjaoli \cite{Me} gave an explicit expression of the $(k,a)$-generalized translation of $f$. From the explicit expression we have
\begin{align*}M_f(x,t)&=\frac12\left(f(x\ast t)+f(x\ast-t)\right)\\&=\frac{M_{k,n}}{2n}\left(\int_0^\pi f_{\mathrm e}\left(\langle\langle x,t\rangle\rangle_{\phi,n}\right)(\sin\phi)^{2nk-n}d\phi+\int_0^\pi f_{\mathrm o}\left(\langle\langle x,t\rangle\rangle_{\phi,n}\right)K(x,t)(\sin\phi)^{2nk-n}d\phi\right),\end{align*}
where $$K(x,t)=\frac{n!\operatorname{sgn}(x)}{(2kn-n)_n}C_n^{nk-\frac n2}\left(\frac{|x|^{\frac{1}{n}}-t^{\frac{1}{n}} \cos \phi}{\langle\langle x, t\rangle\rangle_{\phi, n}^{\frac{1}{n}}}\right),$$
$\langle\langle x,t\rangle\rangle_{\phi,n}:=\left(\vert x\vert^\frac2n+ t^\frac2n-2\vert xt\vert^\frac1n\cos\phi\right)^\frac n2$, $M_{k, n}=\frac{n \frac{n(2 k-1)}{2}}{2^{\frac{n(2 k-1)+2}{2}} \Gamma\left(\frac{n(2 k-1)+2}{2}\right)} $,
and $C_n^\lambda$ are the Gegenbauer polynomials defined by the generating function
$$\left(1-2 u r+r^2\right)^{-2 \lambda}=\sum_{n=0}^{\infty} C_n^\lambda(u) r^n.$$
The Gegenbauer polynomial $C_n^\lambda$ takes the explicit form
$$C_n^\lambda(u)=\frac{1}{\Gamma(\lambda)} \sum_{k=0}^{[n / 2]}(-1)^k \frac{\Gamma(n-k+\lambda)}{k!(n-2 k)!}(2 u)^{n-2 k}.$$
Obviously,
$$\left|\frac{|x|^{\frac{1}{n}}-t^{\frac{1}{n}} \cos \phi}{\langle\langle x, t\rangle\rangle_{\phi, n}^{\frac{1}{n}}}\right|\leq 1$$
Then from the estimate 
$$\left|C_n^\lambda(u)\right|\leq C_n^\lambda(1)=\frac{(2\lambda)_n}{n!}\;\;\;\;\mathrm{for}\;\;-1\leq u\leq1,$$ we have $\left|K(x,t)\right|\leq1$.
Therefore, 
\begin{align*}M_f(x,t)&\geq\frac{M_{k,n}}{2n}\Bigg(\int_0^\pi f_{\mathrm e}\left(\mathrm{sgn}\left(K(x,t)\right)\langle\langle x,t\rangle\rangle_\phi\right)\left|K(x,t)\right|(\sin\phi)^{2kn-n}d\phi\\&\;\;\;\;\;\;\;\;\;\;+\int_0^\pi f_{\mathrm o}\left(\mathrm{sgn}\left(K(x,t)\right)\langle\langle x,t\rangle\rangle_\phi\right)\mathrm{sgn}\left(K(x,t)\right)K(x,t)(\sin\phi)^{2kn-n}d\phi\Bigg)\\&=\frac{M_{k,n}}{2n}\left(\int_0^\pi\left(f_{\mathrm e}+f_{\mathrm o}\right)\left(\pm\langle\langle x,t\rangle\rangle_\phi\right)\left|K(x,t)\right|(\sin\phi)^{2kn-n}d\phi\right),\end{align*}
which is obviously positivity-preserving.

We then consider for $N\geq 2$, the following property for the integral kernel $B_{k,a}(x,y)$ of the $(k,a)$-generalized Fourier transform was given in \cite{GIT3}.

\begin{prop}\label{Heck}
(\cite[Proposition 5.7]{GIT3}) If $x,y\in\mathbb{R}^N$, $x=\rho x'$, $y=vy'$, then
$$\frac{\displaystyle1}{\displaystyle d_{k,a}}\int_{\mathbb S^{N-1}}B_{k,a}(x,vy')p(y'){\vartheta_{k,a}}\left( y'
\right)\,d\sigma(y')=\frac{e^{-\frac{i\pi m}a}\Gamma(\lambda_a+1)}{a^{2m/a}\,\Gamma(\lambda_{k,a,m}+1)}\,v^mj_{\lambda_{k,a,m}}\left(\frac2a\,{(\rho v)}^{a/2}\right)p(x),$$
where $p$ is a polynomial of degree $m$ and $j_\lambda(x)=2^\lambda\Gamma(\lambda+1)x^{-\lambda}J_\lambda(x)$ is the normalized Bessel function and $J_\lambda(x)$ is the classical Bessel
function.

\end{prop}
For $f \in\mathcal L^1_k(\mathbb R^N)$, we have
\begin{align}\label{Mff}M_f(x,t)&=\frac{c_{k,1}}{d_{k,1}}\int_{\mathbb S^{N-1}}\left(\int_{\mathbb R^N}B_{k,1}\left(x,\xi\right)B_{k,1}\left(ty,\xi\right)F_{k,1}\left(f\right)\left(\xi\right)\vartheta_{k,1}\left(\xi\right)d\xi\right)\vartheta_{k,1}(y)d\sigma(y)\notag\\&=\frac{\displaystyle c_{k,1}}{\displaystyle d_{k,1}}\int_{\mathbb R^N}B_{k,1}\left(x,\xi\right)F_{k,1}\left(f\right)\left(\xi\right)\vartheta_{k,1}\left(\xi\right)d\xi\int_{\mathbb S^{N-1}}B_{k,1}\left(ty,\xi\right)\vartheta_{k,1}(y)d\sigma(y)\notag\\&=c_{k,1}\int_{\mathbb R^N}B_{k,1}\left(x,\xi\right)j_{2\left\langle k\right\rangle+N-2}\left(2\sqrt{t\left|\xi\right|}\right)F_{k,1}\left(f\right)\left(\xi\right)\vartheta_{k,1}\left(\xi\right)d\xi,\end{align}
where $B_{k,1}(x,\xi)$ has the expression (see \cite[Theorem 4.24]{BSKO})
$$B_{k,1}(x,\xi)=\Gamma\left(\frac{N-1}2+\langle k\rangle\right)V_k\left[{\widetilde J}_{\frac{N-3}2+\langle k\rangle}\left(\sqrt{2\left|x\right|\left|\xi\right|\left(1+\left\langle\frac x{\left|x\right|},\cdot\right\rangle\right)}\right)\right]\left(\frac\xi{\left|\xi\right|}\right),$$
and ${\widetilde J}_\lambda(x):=\left(\Gamma(\lambda+1)\right)^{-1}j_\lambda(x)$.
It reduces to the definition of the generalized spherical mean in \cite[Section 6]{GIT2} for the one dimensional case.
Thus
$$F_{k,1}\left(M_f(\cdot,t)\right)\left(\xi\right)=j_{2\left\langle k\right\rangle+N-2}\left(2\sqrt{t\left|\xi\right|}\right)F_{k,1}\left(f\right)\left(\xi\right).$$
If $2\langle k\rangle+N-2>0$, it can be observed that $M_f\in C(\mathbb{R}^N\times{\mathbb{R}}_+)$. Moreover, if $f\in \widehat{\mathcal S}(\mathbb{R}^N)$, then $M_f\in C^\infty(\mathbb{R}^N\backslash\{0\}\times{\mathbb{R}}_+)$ and $M_f(\cdot,t)\in L^2\left(\mathbb{R}^N,\vartheta_{k,1}\left(x\right)dx\right)$.

In \cite[Section 4]{BD}, S. Ben Sa\"id and L. Deleaval studied the following heat kernel of the heat operator $\left|x\right|\Delta_k-\partial_t$ with $x=\left|x\right|\theta'$ and $y=\left|y\right|\theta''$,
$$
h_k(x, y ; t)=\frac{c_{k,1}}{t^{N+2\langle k\rangle-1}} \Gamma\left(\frac{N-1}{2}+\langle k\rangle\right) \mathrm{e}^{-\left(\frac{|x|+|y|}{t}\right)} V_k\left[\widetilde{I}_{\frac{N-3}{2}+\langle k\rangle}\left(\frac{\sqrt{2|x||y|\left(1+\left\langle\theta^{\prime}, \cdot\right\rangle\right)}}{t}\right)\right]\left(\theta^{\prime \prime}\right),
$$
where $\widetilde{I}_{\lambda}(\omega)$ is the $I$-normalized Bessel function 
${\widetilde I}_\lambda(\omega):=\left(\Gamma(\lambda+1)\right)^{-1}j_\lambda(i\omega)$.
The kernel $h_k(x, y ; t)$ is strictly positive and it is obvious that $h_k\left(\cdot,y;t\right)\in\mathcal L^1_k(\mathbb R^N)$. 
And they got the following formula from the definitions of the translation operator
\begin{align}\label{heat kernel}F_{k,1}\left(h_k\left(\cdot,y;t\right)\right)(\xi)=c_{k,1}^{}\mathrm e^{-t\left|\xi\right|}B_{k,1}(y,\xi),\end{align}
since $h_k(x, y ; t)$ can be written as
$$
c_{k,1}^2 \int_{\mathbb{R}^N} \mathrm{e}^{-t|\xi|} B_{k,1}(x, \xi) B_{k,1}(y, \xi) \vartheta_{k,1}(\xi) d \xi.
$$
The proof of the positivity-preserving property of $M_f$ reduces to the positivity on the heat kernel from the following lemma.
\begin{lem}\label{heatkernel}
Let $f\in \widehat{\mathcal S}(\mathbb{R}^N)$ and $(x,t)\in\mathbb{R}^N\times{\mathbb{R}}_+$. Then
$$\lim_{s\rightarrow0}\int_{\mathbb{R}^N}M_{h_k(\cdot,z;s)}(x,t)f(z)\vartheta_{k,1}(z)dz=M_f(x,t).$$
\end{lem}

\begin{proof}
Denote $g_{x,y}(\xi):=c_{k,1}^{}\mathrm e^{-s\left|\xi\right|}B_{k,1}(x,\xi)B_{k,1}(y,\xi)\in {{L}^{2}}( {{\mathbb{R}}^{N}},$ ${\vartheta_{k,1}}\left( x
\right)dx ).$
Then from the inversion formula of $(k,1)$-generalized Fourier transform, 
$$h_k\left(x\ast y,z;s\right)=F_{k,1}\left(g_{x,y}\right)(z).$$
For a function $f\in \widehat{\mathcal S}(\mathbb{R}^N)$, $f,\;F_{k,1}(f)\in {(L^1\cap L^2)}( {{\mathbb{R}}^{N}},$ ${\vartheta_{k,1}}\left( x
\right)dx )$.
Thus we have
\begin{align*}\int_{\mathbb{R}^N}M_{h_k(\cdot,z;s)}(x,t)f(z)\vartheta_{k,1}(z)dz&=\frac{c_{k,1}^{}}{d_{k,1}}\int_{\mathbb{R}^N}\int_{\mathbb S^{N-1}}F_{k,1}(g_{x,ty})(z)\vartheta_{k,1}(y)d\sigma(y)f(z)\vartheta_{k,1}(z)dz\\&=\frac{c_{k,1}^{}}{d_{k,1}}\int_{\mathbb S^{N-1}}\int_{\mathbb{R}^N}g_{x,ty}(z)F_{k,1}\left(f\right)(z)\vartheta_{k,1}(z)dz\vartheta_{k,1}(y)d\sigma(y).\end{align*}
Then by dominated convergence theorem,
\begin{align*}\lim_{s\rightarrow0}\int_{\mathbb{R}^N}M_{h_k(\cdot,z;s)}(x,t)f(z)\vartheta_{k,1}(z)dz&=\frac{c_{k,1}^{}}{d_{k,1}}\int_{\mathbb S^{N-1}}\int_{\mathbb{R}^N}B_{k,1}(x,z)B_{k,1}(ty,z)F_{k,1}\left(f\right)(z)\vartheta_{k,1}(y)\vartheta_{k,1}(z)dzd\sigma(y)\\&=\frac1{d_{k,1}}\int_{\mathbb S^{N-1}}f(x\ast ty)\vartheta_{k,1}(y)d\sigma(y)=M_f(x,t).\end{align*}
\end{proof}

The following product formula of the integral kernel of $(k,1)$-generalized translation was obtained in the proof of the radial formula of the $(k,1)$-generalized translation in \cite{BD}.
\begin{prop}(\cite{BD})
For $x=r'\theta',\;z=r''\theta''$,
\begin{align*}
\int_{\mathbb S^{N-1}}B_{k,1}(x,&r\omega)B_{k,1}(z,r\omega)\vartheta_{k,1}(\omega)d\omega=c_{k,1}^{-1}\frac{\Gamma\left(\frac{N-1}2+\langle k\rangle\right)}{\sqrt\pi\Gamma\left(\frac{N-2}2+\langle k\rangle\right)}\times\\&\;\;\;\;V_k\left[\int_{-1}^1j_{N+2\langle k\rangle-2}\left(2\sqrt{r'r+r''r-r\sqrt{2r'r''\left(1+\left\langle\theta',\cdot\right\rangle\right)}u}\right)\left(1-u^2\right)^{\frac N2+\langle k\rangle-2}du\right]\left(\theta''\right).
\end{align*}
\end{prop}

\noindent{\it Proof of Theorem \ref{positivity} 1)}.
Let us consider for a function $f\in C_0(\mathbb{R}^N\backslash\left\{0\right\})$ first. 
It suffices to prove for $f\in W_{k,1}(\mathbb R^N)$ because $W_{k,1}(\mathbb R^N)$ is ${\left\|\cdot\right\|}_\infty$-dense in $C_0(\mathbb{R}^N\backslash{\left\{0\right\}})$. 
It remains to show that $M_{h_k(\cdot,z;s)}(x,t)\geq 0$ according to Lemma \ref{heatkernel}.
Denote $M(x,t):=M_{h_k(\cdot,z;s)}(x,t)$ for brevity. Then from \eqref{Mff}, \eqref{heat kernel} and polar coordinate transformation, 
\begin{align*}M(x,t)&=c_{k,1}^2\int_{\mathbb{R}^N}\mathrm e^{-s\left|\xi\right|}B_{k,1}(x,\xi)B_{k,1}(z,\xi)j_{2\langle k\rangle+N-2}\left(2\,\left|t\xi\right|^\frac12\right)\vartheta_{k,1}(\xi)d\xi\\&\overset{\xi=r\omega}=c_{k,1}^2\int_0^\infty\mathrm e^{-sr}I(x,z,r) j_{2\langle k\rangle+N-2}\left(2\,\left|tr\right|^\frac12\right)r^{2\left\langle k\right\rangle+N-2}dr,\end{align*}
where
\begin{align*}I(x,z,r)&=\int_{\mathbb S^{N-1}}B_{k,1}(x,r\omega)B_{k,1}(z,r\omega)\vartheta_{k}(\omega)d\omega\\&=c_{k,1}^{-1}\frac{\Gamma\left(\frac{N-1}2+\langle k\rangle\right)}{\sqrt\pi\Gamma\left(\frac{N-2}2+\langle k\rangle\right)}\times\\&\;\;\;\int_{\mathbb{R}^{\mathrm N}}\int_{-1}^1j_{N+2\langle k\rangle-2}\left(2\sqrt{r\left(\left|x\right|+\left|z\right|-\sqrt{2\left(\left|x\right|\left|z\right|+\left\langle x,\eta\right\rangle\right)}u\right)}\right)\left(1-u^2\right)^{\frac N2+\langle k\rangle-2}dud\mu_z(\eta).
\end{align*}
Put $$v_{z,u}(\eta):=\left|x\right|+\left|z\right|-\sqrt{2\left(\left|x\right|\left|z\right|+\left\langle x,\eta\right\rangle\right)}u.$$
Then from the well-known product formula (see e.g. \cite[3.5.61]{BH}) for Bessel functions $j_\alpha$ with $\alpha> -1/2$
\begin{equation*}
j_\alpha(u w) j_\alpha(v w)=\int_0^{\infty} j_\alpha(\xi w) d \nu_{u, v}^\alpha(\xi) \quad \text { for all } w \in \mathbb C,
\end{equation*}
where $d\nu_{u, v}^{\alpha}$ are probability measures on $\mathbb R_+$,
we have
\begin{align*}M(x,t)&=c_{k,1}^{}\frac{\Gamma\left(\frac{N-1}2+\langle k\rangle\right)}{\sqrt\pi\Gamma\left(\frac{N-2}2+\langle k\rangle\right)}\int_0^\infty\mathrm e^{-sr}\\&\;\;\;\int_{\mathbb{R}^{\mathrm N}}\int_{-1}^1j_{2\langle k\rangle+N-2}\left(2\left(rv_{z,u}(\eta)\right)^\frac12\right)j_{2\langle k\rangle+N-2}\left(2\,\left|tr\right|^\frac12\right)\left(1-u^2\right)^{\frac N2+\langle k\rangle-2}r^{2\left\langle k\right\rangle+N-2}dud\mu_z(\eta)dr\\&=c_{k,1}^{}\frac{\Gamma\left(\frac{N-1}2+\langle k\rangle\right)}{\sqrt\pi\Gamma\left(\frac{N-2}2+\langle k\rangle\right)}\int_0^\infty e^{-sr}\\&\;\;\;\int_{\mathbb R^N}\int_{-1}^1\left(\int_0^\infty j_{2\langle k\rangle+N-2}\left(2\xi r^\frac12\right)d\nu_{\sqrt{v_{z,u}{(\eta)}},\sqrt t}^{2\langle k\rangle+N-2}(\xi)\right)\left(1-u^2\right)^{\frac N2+\langle k\rangle-2}r^{2\left\langle k\right\rangle+N-2}dud\mu_z(\eta)dr.\end{align*}
By a change of variable in the formula 11.4.29 of \cite{AS}, we get
$$\frac1{2\Gamma\left(v+1\right)}\int_0^\infty e^{-at}t^vj_v(b\sqrt t)dt=\frac1{a^{v+1}}e^{-\frac{b^2}{4a}},\;\;\mathfrak Rv>-1,\;\mathfrak Ra>0.$$
Therefore,
\begin{align*}M(x,t)=d_{k,1}^{}\frac{2\Gamma\left(\frac{N-1}2+\langle k\rangle\right)}{\sqrt\pi\Gamma\left(\frac{N-2}2+\langle k\rangle\right)}\int_{\mathbb R^N}\int_{-1}^1\left(\int_0^\infty\frac1{s^{v+1}}e^{-\frac{\xi^2}s}d\nu_{\sqrt{v_{z,u}{(\eta)}},\sqrt t}^{2\langle k\rangle+N-2}(\xi)\right)\left(1-u^2\right)^{\frac N2+\langle k\rangle-2}dud\mu_z(\eta).\end{align*}
This integral is non-negative obviously and the operator $f\mapsto M_f$ is therefore positivity-preserving on $C_0(\mathbb{R}^N\backslash\left\{0\right\})$.

For a non-negative function $f\in C_0(\mathbb{R}^N)$, if $f(0)=0$, then the proof reduces to that for functions in $C_0(\mathbb{R}^N\backslash\left\{0\right\})$.
If $f(0)>0$, we define $\widetilde{f}(r):=\underset{x'\in \mathbb S^{N-1}}{\min}f(rx')$ for $x=rx'$ and a radial function $\widetilde{f}_0(x):=\widetilde{f}(\left|x\right|)$. Then $0\leq \widetilde{f}_0(x)\leq f(x)$ and $\widetilde{f}_0(0)=f(0)$. Let us write $f=(f-\widetilde{f}_0)+\widetilde{f}_0$. It can then be observed that the spherical mean operator is positivity-preserving on $f$ from its positivity-preserving property on $C_0(\mathbb{R}^N\backslash\left\{0\right\})$ along with the positivity-preserving property of the generalized translation operator on radial functions.$\hfill\Box$

\section{The representing measure $\sigma_{x,t}^{k,1}$ of the generalized spherical mean}

For each $x, y \in \mathbb{R}^N$, consider the linear functional $\Psi_{x,y}:f\mapsto f\left(x\ast y\right)$ on $\left(C_{0, \operatorname{rad}}\left(\mathbb{R}^N\right),\|\cdot\|_{\infty}\right)$. It is positive and bounded with norm $\left\|\Psi_{x, y}\right\|=1$ according to the formula \eqref{radial} for radial functions. 
We can then easily obtain the following proposition parallel to the Theorem 5.1 in \cite {R2}.
\begin{prop}
 For each $x, y \in \mathbb{R}^N$ there exists a unique compactly supported, radial probability measure $\rho_{x, y}^{k,1}\in M^1(\mathbb R^N)$ such that for all $f \in C_{\text {rad }}\left(\mathbb{R}^N\right)$,

$$
f\left(x *_k y\right)=\int_{\mathbb{R}^N} f d \rho_{x, y}^{k,1}.
$$
The support of $\rho_{x, y}^{k,1}$ is contained in

$$
\left\{\xi\in\mathbb{R}^N:d_G{(x,y)}^2\leq\left|\xi\right|\leq\underset{g\in G}{max}\left(\left|x\right|+\left|y\right|+\sqrt{2\left(\left|x\right|\left|y\right|+\left\langle x,y\right\rangle\right)}\right)\right\},
$$
where $d_G(x,y)$ is the distance of orbits $d_G(x,y):
=\underset{g\in G}{\min}\;d\left(gx,y\right)$.
\end{prop}

We can also get the following proposition according to the radial formula \eqref{radial}.
Below we denote $B\left(x,r\right):=\left\{y:d\left(y,x\right)\leq r\right\}$ and $B_r:=\left\{y\in\mathbb{R}^N:\sqrt{\left|y\right|}\leq
r\right\}$.
\begin{prop}\label{suppf}
Let $f=f_0\left(\left|\cdot\right|\right)$ be a radial function on ${{L}^{2}}\left(
{{\mathbb{R}}^{N}},{\vartheta_{k,1}}\left( x \right)dx \right)$,
$\operatorname{supp} f\subseteq B_r$, then $$\operatorname{supp}\tau_xf\subseteq \bigcup_{g\in
G}B(gx,r).$$
\end{prop}

\begin{proof}
It is easy to see that for any $\eta\in
co(G.x)$ and $u\in[-1,1]$, we have
\begin{align}\label{min}\sqrt{\left|x\right|+\left|y\right|-\sqrt{2\left(\left|x\right|\left|y\right|+\left\langle\eta,y\right\rangle\right)}u}\geq \underset{g\in G}{\min}\;d\left(gx,y\right).\end{align}
Then we can derive the proposition by firstly observing for for all continuous radial functions and then using a density argument and the continuity of of the $(k,1)$-generalized translation on ${{L}^{2}}\left(
{{\mathbb{R}}^{N}},{\vartheta_{k,1}}\left( x \right)dx \right)$.
\end{proof}

We then investigate the  measure representing the generalized spherical mean $M_f(x,t)$.
For fixed $x\in \mathbb R^N$ and $t\geq 0$, consider the linear functional
\[\Phi_{x,t}: \,f\mapsto M_f(x,t).\]
It is a positive linear functional for $a=1$ according to the first part of Theorem \ref{positivity}.
Moreover, $\Phi_{x,t}(1)=1$.
It follows that the linear functional $\Phi_{x,t}$ is represented by
a compactly supported  probability  measure $\sigma_{x,t}^{k,1}\in M^1(\mathbb R^N)$ (cf. \cite[Theorem 2.1.7]{Ho}), i.e.,

\begin{align} \label{measure}M_f(x,t) = \int_{\mathbb R^N} f\,d\sigma_{x,t}^{k,1} \quad\text{for all } \,f\in
  C_0(\mathbb{R}^N).\end{align}

\noindent{\it Proof of Theorem \ref{positivity} 2)}.  
The transformation properties \eqref{transformation} of $\sigma_{x, t}^{k,1}$ can be deduced immediately from the  invariance properties \eqref{kernelproperty} of the integral kernel $B_{k,a}(x,y)$.
We will then analyze the support of $\sigma_{x, t}^{k,1}$.

For the one dimensional case, we have
$$\sigma_{x,t}^{k,1}=\frac12\left(\nu_{x,t}^{k,1}+\nu_{x,-t}^{k,1}\right)$$
according to \eqref{tauN=1}.
Obviously, $\operatorname{supp}\sigma_{x,t}^{k,1}\subseteq\left\{\xi\in\mathbb{R}:\vert\sqrt{\vert x\vert}-\sqrt t\vert\leq\sqrt{\vert\xi}\vert\leq\sqrt{\vert x\vert}+\sqrt t\right\}.$

For $N\geq 2$, we
firstly show that   
 the measures $\sigma_{x, t}^{k,1}$ satisfy

\begin{align}\label{suppf1}
\operatorname{supp} \sigma_{x, t}^{k,1} \subseteq\left\{\xi \in \mathbb{R}^N:\sqrt{|\xi|} \geq \vert\sqrt{\vert x\vert}-\sqrt t\vert\right\}.
\end{align}
Suppose to the contrary that $\operatorname{supp}\sigma_{x,t}^{k,1}\not\subseteq\left\{\xi\in\mathbb{R}^N:\vert\xi\vert\geq\vert\sqrt{\vert x\vert}-\sqrt t\vert^2\right\}$. Then there exists some radial function $f \in   C_0(\mathbb{R}^N)$ with $f \geq 0$, 

\begin{align}\label{measuresupport}
\operatorname{supp} f \cap\left\{\xi \in \mathbb{R}^N:|\xi| \geq \vert\sqrt{\vert x\vert}-\sqrt t\vert^2\right\}=\varnothing
\end{align}
and such that $M_f(x, t)>0$. But then $\eta \mapsto f\left(x * t \eta\right)$ is not identically zero on $\mathbb S^{N-1}$. Then in view of Proposition 5.1, we have
$$\operatorname{supp}f\cap\left\{d_G{(x,t\eta)}^2\leq\left|\xi\right|\leq\underset{g\in G}{max}\left(\left|x\right|+t+\sqrt{2\left(\left|x\right|+\left\langle x,g\eta\right\rangle\right)t}\right)\right\}\neq\varnothing.$$
But this contradicts with \eqref{measuresupport}.

Then, we will show that 
$$\operatorname{supp}\sigma_{x,t}^{k,1}\subseteq\bigcup_{g\in G}B(gx,\sqrt t).$$
The parallel result in \cite{R2} uses techniques on the energy method in partial differential equations by studying the domain of dependence of a wave equation first. However, such energy method cannot be applied to our case because the solution to the relevant deformed wave equation (see \eqref{wave}) may not be differentiable at least at the origin due to its non-smooth coefficient $\left|x\right|$, and the derivative of the energy integral does not make sense due to the singularity at the origin. But we can find another way of proof here. Similar method also applies to the proof of Theorem 4.1 in \cite{R2}.

From the polar coordinates transformation $y=sx'$, $0\leq s\leq t$, $x'\in \mathbb S^{N-1}$,
$$\int_{B_{\sqrt t}}\tau_xf(y)\vartheta_{k,1}(y)dy=\int_0^t\int_{\mathbb S^{N-1}}\tau_xf(sx')s^{2\left\langle k\right\rangle+N-2}\vartheta_{k,1}(x')d\sigma(x')ds=d_{k,1}\int_0^ts^{2\left\langle k\right\rangle+N-2}M_f(x,s)ds$$
If $f$ is a nonnegative function in $\mathcal L^1_k(\mathbb R^N)$ and
$$\text{supp}\,f\cap\left(\bigcup_{g\in G}B(gx,\sqrt t)\right)=\varnothing,$$
then from the symmetry property of the $(k,1)$-generalized translation and Proposition \ref{suppf},
$$\int_{\mathbb{R}^N}1_{B_{\sqrt t}}(y)\tau_xf(y)\vartheta_{k,1}(y)dy=\int_{\mathbb{R}^N}\tau_x1_{B_{\sqrt t}}(y)f(y)\vartheta_{k,1}(y)dy=0.$$ Thus $M_f(x,s)=0$ for all $0\leq s\leq t$, since $M_f(x,s)\geq 0$ and $M_f\in C(\mathbb{R}^N\times{\mathbb{R}}_+)$.
And so $$\int_{\mathbb{R}^N}f\,d\sigma_{x,t}^{k,1}\,=\,M_f(x,t)=0.$$
This completes the proof of Theorem \ref{positivity} 2).
  $\hfill\Box$

\noindent{\it Proof of Theorem \ref{positivity} 3)}. 
The definition of the $(k,1)$-generalized spherical mean operator can then be extended to all function in $C_b(\mathbb{R}^N)$ since the representing measure $\sigma_{x,t}^{k,1}$ is compactly supported. And if we take $f=B_{k,1}(\cdot,\xi)$ in \eqref{spherical mean},  then from Proposition \ref{Heck} we have
\begin{align}\label{MfB}M_f(x,t)=B_{k,1}(x,z)j_{2\left\langle k\right\rangle+N-2}\left(2\sqrt{t\left|z\right|}\right).\end{align}
From \eqref{MfB}, \eqref{Mff}, \eqref{Borel measure} and \eqref{inversion}, \eqref{measure} is equivalent to the following product formula
$$B_{k,1}(x,z)j_{2\left\langle k\right\rangle+N-2}\left(2\sqrt{t\left|z\right|}\right)=\int_{\mathbb R^N} B_{k,1}(\xi,z)\,d\sigma_{x,t}^{k,1}(\xi).$$
Thus the measure $\sigma_{x,t}^{k,1}$ representing the spherical mean is unique from the injectivity of the $(k,1)$-generalized Fourier transform on bounded Borel measures. 
  $\hfill\Box$

\

\noindent{\it Proof of Theorem \ref{wave equation}}

The uniqueness of the solution to the equation \eqref{wave} for $f\in C_c^\infty(\mathbb{R}^N)$ was proven in \cite{Ba}. However, it is not a rigorous proof  in that it relied on a false result in \cite{Jo} that  the Schwartz space is invariant under $F_{k,1}$. And we also lack results on properties of the $(k,1)$-generalized translation of compactly supported smooth functions. But the proof of Lemma 1 in \cite{Ba} makes sense once we make a modification. To ensure that the solution $\mathcal R_{\lambda_1}^{-1}M_f(\cdot,\frac12t^2)\in L^2\left(\mathbb{R}^N,\vartheta_{k,1}\left(x\right)dx\right)$, where $\lambda_1=2\left\langle k\right\rangle+N-2$ and $\mathcal R_{\lambda_1}^{-1}M_f\in C(\mathbb{R}^N\times\mathbb R_+)\cap C^\infty(\mathbb{R}^N\backslash\{0\}\times\mathbb R_+)$, we modify the space of the initial function $f$ to be $\widehat{\mathcal S}(\mathbb{R}^N)$ rather than the  compactly supported smooth function space.

To show that $\mathcal R_{\lambda_1}^{-1}M_f(x,\frac12t^2)$ is the solution, we consider the
singular Sturm-Liouville operator
 $A_\alpha^t$ for $\alpha\geq -1/2,$
\[ A_\alpha^t\,:=\, \partial_t^2 + \frac{2\alpha+1}{t}\partial_t,\;t>0.\]
For fixed $z\in\mathbb C^N$, the Bessel functions $j_\alpha(tz)$ are eigenfunctions of the Sturm-Liouville operator (see \cite{R2}).
By substituting $t$ by $\sqrt{2t}$ ($\sqrt{\frac2a}t^\frac a2$ for $a=1$), we get the deformed Sturm-Liouville operator
\begin{align*}A_{1,\alpha}^t:&=\left(2t\partial_t^2+\partial_t\right)+(2\alpha+1)\partial_t\;\\&=2\left(t\partial_t^2+(\alpha+1)\partial_t\right).\end{align*}
And for fixed $z\in\mathbb C^N$, the Bessel functions $j_\alpha\left(2\,\sqrt{t\left|z\right|}\right)$ are  eigenfunctions of the deformed Sturm-Liouville operator, i.e.,
$$A_{1,\alpha}^tj_\alpha\left(2\,\sqrt{t\left|z\right|}\right)=\,-2\left|z\right|j_\alpha\left(2\,\sqrt{t\left|z\right|}\right).$$
Combining with \eqref{Mff} and the fact that $B_{k,1}\left(x,\cdot\right)$ is eigenfunction of the operator $\left|x\right|\triangle_k^x$,  we conclude $u=M_f(x,t)$, $f\in \widehat{\mathcal S}(\mathbb{R}^N)$ is a solution of the equation for the Darboux-type differential-reflection operator $2\left|x\right|\Delta_k^x -A_{1,\lambda_1}^t$, i.e.,
\begin{align} \label{(1.8)}
 &(2\left|x\right|\Delta_k^x - A_{1,\lambda_1}^t) u \,=\, 0 \quad\text{in }\,
\mathbb{R}^N\backslash\{0\}\times \mathbb R_+\,;\\
 &u(x,0) = f(x), \,\, u_t (x,0) = 0 \quad\text{for all } x\in \mathbb R^N.\notag
\end{align}

We then involve the Riemann-Liouville
transform with parameter $\alpha >-1/2$ on $\mathbb R_+$. It is given by
\begin{equation}\label{(1.8a)}
 \mathcal R_\alpha f(t)\,=\,
\frac{2\Gamma(\alpha+1)}{\Gamma(1/2)\Gamma(\alpha +1/2)}
\int_0^1 f(st) (1-s^2)^{\alpha-1/2}ds
\end{equation}
for $f\in C^\infty(\mathbb R_+),$ see \cite{T}.
The operator $\mathcal R_\alpha$  satisfies the intertwining property
\begin{equation}\label{intert1} A_\alpha \mathcal R_\alpha\,=\,\mathcal R_\alpha\, \frac{d^2}{dt^2}.\end{equation}
Substituting $t$ by $s=\sqrt{2t}$ in \eqref{intert1}, we get
\begin{equation}\label{intert}A_{1,\alpha}\mathcal R_\alpha\,=\mathcal R_\alpha\frac{\displaystyle d^2}{\displaystyle ds^2}\,=\,\mathcal R_\alpha\left(2t\frac{\displaystyle d^2}{\displaystyle dt^2}+\frac{\displaystyle d}{\displaystyle dt}\,\right).\end{equation}

Put
 $ u_f(x,t):=
\mathcal R_{\lambda_1}^{-1}M_f(x,t)$. Then according to \eqref{(1.8)} and the intertwining property 
\eqref{intert}, $ u=
  u_f$  solves the initial value problem
\begin{align}\label{(1.40)}
 &(2\left|x\right|\triangle_k  -  \left(2t\partial_t^2+\partial_t\,\right)) u \,=\, 0 \quad\text{in }\, \mathbb R^N\backslash\{0\}\times
\mathbb R_+;\notag\\
 &u(x,0) = f(x), \,\, \sqrt{2t}u_t (x,0) = 0 \quad\text{for all } x\in \mathbb R^N.
\end{align}
Substituting $t$ by $\frac12t^2$ in \eqref{(1.40)}, we get 
$\mathcal R_{\lambda_1}^{-1}M_f(x,\frac12t^2)$ is the solution to 
\begin{align}\label{(1.41)}
 &u_{tt}-2\left|x\right|\triangle_ku=0\;\;(x,t)\in\mathbb{R}^N\backslash\{0\}\times
\mathbb R_+;\notag\\
 &u(x,0) = f(x), \,\, u_t (x,0) = 0 \quad\text{for} \;\; f\in
\widehat{\mathcal S}(\mathbb{R}^N), \;\;x\in \mathbb R^N.
\end{align}
$\hfill\Box$

\noindent{\it Proof of Theorem \ref{domain of dependence}}

For fixed $(x,t)\in\mathbb{R}^N\backslash\{0\}\times
\mathbb R_+$,
suppose $f$ is $0$ in $\bigcup_{g\in G}B(gx,t/\sqrt 2)$. Then $f$ vanishes in $\bigcup_{g\in G}B(gx,\sqrt s)$ for all $0\leq \sqrt s\leq t/\sqrt 2$, and so
$M_f(x,s)=0$, $0\leq s\leq \frac12t^2$.
Therefore, the value of the solution
$\mathcal R_{\lambda_1}^{-1}M_f(x,\frac12t^2)=0$ according to the explicit expression of $\mathcal R_\alpha^{-1}$ in \cite{T}.
$\hfill\Box$

\section{Solution to the deformed wave equation}

Consider the deformed wave equaton
\begin{align*}
 &u_{tt}-2\left|x\right|\triangle_ku=0\;\;(x,t)\in\mathbb{R}^N\backslash\{0\}\times
\mathbb R;\notag\\
 &u(x,0) = f(x), \,\, u_t (x,0) = g(x) \quad\text{for} \;\; f,\;g\in
\widehat{\mathcal S}(\mathbb{R}^N), \;\;x\in \mathbb R^N.
\end{align*}
It has the unique solution $u(\cdot,t)\in C(\mathbb{R}^N)\cap C^\infty(\mathbb{R}^N\backslash\{0\})\cap L^2\left(\mathbb{R}^N,\vartheta_{k,1}\left(x\right)dx\right)$. And  it was found in \cite{Ba} that the solution has the form
\begin{align}\label{sol}u(x,t)={\left\langle P_{k,t}^{11},\tau_xf\right\rangle}_{k,1}+{\left\langle P_{k,t}^{12},\tau_xg\right\rangle}_{k,1},\end{align}
where 
$$P_{k,t}^{11}=F_{k,1}(\cos(t\sqrt{2\left|\cdot\right|})),\;\;\; P_{k,t}^{12}=F_{k,1}(\sin(t\sqrt{2\left|\cdot\right|})/\sqrt{2\left|\cdot\right|})$$
and ${\left\langle f,g\right\rangle}_{k,1}:=\int_{\mathbb{R}^N}f(x)g(x)\vartheta_{k,1}(x)dx.$

To interpret the above form of the solution, we extend the $(k,1)$-generalized Fourier transform $F_{k,1}$ to the distributional sense for $T\in \widehat{\mathcal S}(\mathbb{R}^N)'$ as follows 
$$\left\langle F_{k,1}(T),f\right\rangle_{k,1}:=\left\langle T,F_{k,1}(f)\right\rangle_{k,1}\;\mathrm{for}\;\mathrm{all}\;f\in \widehat{\mathcal S}(\mathbb{R}^N).$$
Here $ \widehat{\mathcal S}(\mathbb{R}^N)'$ stands for the dual of $\widehat{\mathcal S}(\mathbb{R}^N)$.
Then we define the convolution of a distribution $T\in \widehat{\mathcal S}(\mathbb{R}^N)'$ and a function $f\in \widehat{\mathcal S}(\mathbb{R}^N)$ as 
$$T\ast f:={\left\langle T,\tau_xf\right\rangle}_{k,1}.$$
And from the fact that $F_{k,1}^2=\operatorname{id}$ (\cite[Corollary 5.2]{BSKO}) and the definition of the $(k,1)$-generalized translation \eqref{deftranslation}, we have 
$$\left|{\left\langle T,\tau_xf\right\rangle}_{k,1}\right|=\left|{\left\langle F_{k,1}(T),F_{k,1}(\tau_xf)\right\rangle}_{k,1}\right|\leq{\left\langle\left|F_{k,1}(T)\right|,\left|F_{k,1}(f)\right|\right\rangle}_{k,1}$$
since $\left|B_{k,1}(x,y)\right|\leq 1$ if $2\left\langle k\right\rangle+N-2>0$. Thus, the definition $T\ast f$ makes sense for all 
$T\in \widehat{\mathcal S}(\mathbb{R}^N)'$ and $f\in \widehat{\mathcal S}(\mathbb{R}^N)$.

We will write the solution \eqref{sol} in terms of the spherical 
mean operator, parallel to the Theorem 3.18 in \cite{BO}.
From polar coordinates transformation and the Bochner type identity for the $(k,a)$-generalized Fourier transform (see \cite[Theorem 5.21]{BSKO}), we have
\begin{align*}u(x,t)&=\int_0^{+\infty}r^{2\left\langle k\right\rangle+N-2}dr\int_{\mathbb S^{N-1}}P_{k,t}^{11}(ry')\tau_xf(ry')\vartheta_k(y')d\sigma(y')\\&\;\;\;\;+\int_0^{+\infty}r^{2\left\langle k\right\rangle+N-2}dr\int_{\mathbb S^{N-1}}P_{k,t}^{12}(ry')\tau_xg(ry')\vartheta_k(y')d\sigma(y')\\&=d_{k,1}\int_0^{+\infty}r^{2\left\langle k\right\rangle+N-2}H_{1,\lambda_1}(\cos(t\sqrt{2\left(\cdot\right)})(r)M_f(x,r)dr\\&\;\;\;\;+d_{k,1}\int_0^{+\infty}r^{2\left\langle k\right\rangle+N-2}H_{1,\lambda_1}(\sin(t\sqrt{2(\cdot)})/\sqrt{2(\cdot)})(r)M_g(x,r)dr,
\end{align*}
where $H_{1,\lambda_1}$ is the Hankel transform 
$$H_{1,\lambda_1}(f)(r)=\frac1{\Gamma(2\left\langle k\right\rangle+N-1)}\int_0^\infty f(s)j_{2\left\langle k\right\rangle+N-2}\left(2\sqrt{rs}\right)s^{2\left\langle k\right\rangle+N-2}\,ds.$$
On the other hand, by a substitution of variables $r$ by $\sqrt{2r}$ and $s$ by $\sqrt{2s}$ in the proof of Theorem 3.18 in \cite {BO} (cf. \cite[p. 331, formula 5, formula 8)]{Er}), we have
\begin{align*}H_{1,\lambda_1}(\cos(t\sqrt{2\left(\cdot\right)})(r)&=\frac1{\Gamma(2\left\langle k\right\rangle+N-1)}\int_0^\infty\cos(t\sqrt{2s})j_{2\left\langle k\right\rangle+N-2}\left(2\sqrt{rs}\right)s^{2\left\langle k\right\rangle+N-2}\,ds\\&=\left\{\begin{array}{lc}\displaystyle{\frac{2\sqrt\pi}{\Gamma(2\left\langle k\right\rangle+N-1)}t\frac{\left(t^2-2r\right)^{-2\left\langle k\right\rangle-N+\frac12}}{\Gamma\left(-2\left\langle k\right\rangle-N+\frac32\right)}}&\mathrm{if}\;0<\sqrt{2r}<t,\\0&\mathrm{if}\;0<t<\sqrt{2r},\end{array}\right.
\end{align*}
and
\begin{align*}H_{1,\lambda_1}(\sin(t\sqrt{2(\cdot)})/\sqrt{2(\cdot)})(r)&=\frac1{\Gamma(2\left\langle k\right\rangle+N-1)}\int_0^\infty\frac{\sin(t\sqrt{2s})}{\sqrt{2s}}j_{2\left\langle k\right\rangle+N-2}\left(2\sqrt{rs}\right)s^{2\left\langle k\right\rangle+N-2}\,ds\\&=\left\{\begin{array}{lc}\displaystyle{\frac{\sqrt\pi}{\Gamma(2\left\langle k\right\rangle+N-1)}\frac{\left(t^2-2r\right)^{-2\left\langle k\right\rangle-N+\frac32}}{\Gamma\left(-2\left\langle k\right\rangle-N+\frac52\right)}}&\mathrm{if}\;0<\sqrt{2r}<t,\\0&\mathrm{if}\;0<t<\sqrt{2r}.\end{array}\right.\end{align*}
If we introduce the Riemann-Liouville distribution $\mathbb S_{\lambda}(x)$ (see \cite{BO}) for $\lambda\in \mathbb C$ and $x\in \mathbb R$, we have
\begin{align*}H_{1,\lambda_1}(\cos(t\sqrt{2\left(\cdot\right)})(r)&=\frac{2\sqrt\pi}{\Gamma(2\left\langle k\right\rangle+N-1)}t\mathbb S_{-2\left\langle k\right\rangle-N+\frac32}(t^2-2r)\\&=\frac{\sqrt\pi}{\Gamma(2\left\langle k\right\rangle+N-1)}\frac d{dt}\left(\mathbb S_{-2\left\langle k\right\rangle-N+\frac52}(t^2-2r)\right)\end{align*}
and
$$H_{1,\lambda_1}(\sin(t\sqrt{2(\cdot)})/\sqrt{2(\cdot)})(r)=\frac{\sqrt\pi}{\Gamma(2\left\langle k\right\rangle+N-1)}\mathbb S_{-2\left\langle k\right\rangle-N+\frac52}\left(t^2-2r\right).$$
We get the following expression of the solution $u(x,t)$.
\begin{thm}
For $2\left\langle k\right\rangle+N-2>0$ and $(x,t)\in\mathbb{R}^N\backslash\{0\}\times
\mathbb R_+$, 
\begin{align*}
u(x,t)&=d_{k,1}\frac{\sqrt\pi}{\Gamma(2\left\langle k\right\rangle+N-1)}\int_0^{{\textstyle\frac12}t^2}r^{2\left\langle k\right\rangle+N-2}\operatorname{sgn}(t)\frac d{dt}\left(\mathbb S_{-2\left\langle k\right\rangle-N+\frac52}(t^2-2r)\right)M_f(x,r)dr\\&\;\;\;\;+\operatorname{sgn}(t)d_{k,1}\frac{\sqrt\pi}{\Gamma(2\left\langle k\right\rangle+N-1)}\int_0^{\textstyle\frac12t^2}r^{2\left\langle k\right\rangle+N-2}\mathbb S_{-2\left\langle k\right\rangle-N+\frac52}\left(t^2-2r\right)M_g(x,r)dr.
\end{align*}
\end{thm}

\section{Appendix: The metric $d(x,y)$}

The metric
$d\left(x,y\right)=\sqrt{\left|x\right|+\left|y\right|-\sqrt{2\left(\left|x\right|\left|y\right|+\left\langle x,y\right\rangle\right)}}$ was shown in \cite{Te3} to be the metric corresponding to the setting of $(k,1)$-generalized Fourier analysis, since the support of the $(k,1)$-generalized translation in Proposition \ref{suppf} is precise if the multiplicity function $k>0$, i.e.,  if $f$ is a nonnegative
radial function in $L^2\left(\mathbb{R}^N,\vartheta_{k,1}\left(x\right)dx\right)$ with $\operatorname{supp} f= B_r$, then 
$$\operatorname{supp}\tau_xf= \bigcup_{g\in
G}B(gx,r).$$

It can be verified that the metric $\sqrt 2 d(x,x_0)$ solves the following equation from basic calculation
\begin{align}\label{metricproperty}\left|\nabla q(x)\right|= 1/\sqrt{2\left|x\right|},\;q(x_0)=0.\end{align}
It also well known that the Riemannian distance determined by the coefficients $1/\sqrt{2\left|x\right|}$, to be denoted by $L(x,x_0)$, is the maximal subsolution to
the equation \eqref{metricproperty}.
Therefore, 
$\sqrt 2 d(x,x_0)\leq L(x,x_0)$.
If $N=1$, $$d(x,y)=\left\{\begin{array}{l}\sqrt{\left|x-y\right|},\;xy\leq0,\\\left|\sqrt{\left|x\right|}-\sqrt{\left|y\right|}\right|,\;xy>0.\end{array}\right.$$
Obviously, it is  not a Riemannian distance because there exists no continuous rectifiable curve between two distinct points (see \cite[Remark 3.2]{Te3}), and $\sqrt 2 d(x,x_0)$ is not the same metric with $L(x,x_0)$ for $N=1$. However, for $N\geq 2$, the situation is different because the cone $\mathbb{R}^N\backslash\{0\}$  is connected  if $N\geq 2$, and the geodesics are the projections of those on the cone
$C_{ +}:=\left\{\left(\zeta_1, \cdots, \zeta_{N+1}\right) \in \mathbb{R}^{N+1}: \zeta_{N+1}>0, \quad \zeta_1^2+\cdots+\zeta_N^2=\zeta_{N+1}^2\right\}$
onto the plane $\mathbb R^N$.

\begin{lem}
For $N\geq 2$, $\sqrt 2 d(x,x_0)$ equals to the Riemannian distance determined by the coefficients $1/\sqrt{2\left|x\right|}$, denoted by $L(x,x_0)$.
\end{lem}

\begin{proof}

We consider the proof for the two dimensional case first. For $x=(x_1,x_2)\in \mathbb R^2$, from the polar coordinates transformation $x_1=r\cos\theta,\;x_2=r\sin\theta$, $0\leq\theta\leq2\mathrm\pi$,
the Riemannian metric
$\frac1{2\left|x\right|}(dx_1^2+dx_2^2)=\frac1{2r}(dr^2+r^2d\theta^2)$.
Let
$u=\sqrt{2r}\cos\frac\theta2,\;v=\sqrt{2r}\sin\frac\theta2$.
Then $du^2+dv^2=\frac1{2r}(dx_1^2+dx_2^2)$. Therefore,
$L(x,x_0)$ equals to the Euclidean distance of the image of the two points $x$ and $x_0$ in the coordinate system $(u,v)$. We denote $x_0=(r_0\cos\theta_0,r_0\sin\theta_0)$, $0\leq\theta_0\leq2\mathrm\pi$. Then 
\begin{align*}L(x,x_0)&=\sqrt{\left(\sqrt{2r}\cos\frac\theta2-\sqrt{2r_0}\cos\frac{\theta_0}2\right)^2+\left(\sqrt{2r}\sin\frac\theta2-\sqrt{2r_0}\sin\frac{\theta_0}2\right)^2}\\&=\sqrt{2\left(r+r_0-2\sqrt{rr_0}\cos\frac{\theta-\theta_0}2\right)}=\sqrt{2\left(\left|x\right|+\left|x_0\right|-2\sqrt{\left|x\right|\left|x_0\right|}\cos\frac\alpha2\right)},\end{align*}
where $\alpha=\arccos\frac{\left\langle
x,x_0\right\rangle}{\left|x\right|\left|x_0\right|}$. 

For higher dimensional cases, since $\sqrt 2 d(x,x_0)$ is the length of a curve between $x$ and $x_0$ in a two dimensional subspace, $L(x,x_0)$ is  also no more than $\sqrt 2 d(x,x_0)$ from its definition.
\end{proof}

\begin{rem}

We can also get the domain of dependence in Theorem \ref{domain of dependence} of the solution to the deformed wave equation \eqref{wave} at $(x_0,t_0)$ for $G$-invariant functions using energy method, when $x_0$ is far away from the origin and $t_0$ is small, i.e., if $\sqrt{2|x_0|}>t_0$, by constructing the energy
$$e\left(t\right)=\frac12\int_{C_t}\left(u_t^2+2\left|x\right|\left|\nabla_xu\right|^2\right)\left|2x\right|^{-1}\vartheta_{k}(x)dx,$$
where $C_t=\{x\in\mathbb{R}^N\,:L(x,x_0)\leq t_0-t\}$.
The proof is similar to that of 
Lemma 4.3 in \cite{R2}.
\end{rem}

\section*{Acknowledgments}
The work was done during the author's postdoctoral research at the University of Tokyo supported by Japan Society for the Promotion of Science. The author would like to thank his supervisor Toshiyuki Kobayashi and also the referee for valuable comments.



\begin{thebibliography}{9}

\bibitem{AS}
M. Abramowitz and I. A. Stegun, Pocketbook of Mathematical Functions. Verlag Harri
Deutsch, Frankfurt/Main, 1984.

\bibitem{B} S. Ben Sa\"id, On the integrability of a representation of $sl(2,\mathbb R)$,  J. Funct. Anal.,  250(2007), 249--264

\bibitem{Be} S. Ben Sa\"id, A product formula and a convolution structure for a $k$-Hankel transform on R. J. Math.
Anal. Appl., 463(2)(2018), 1132--1146.

\bibitem{Ba}
 S. Ben Said, S. al-Blooshi, M. al-Kaabi, A. al-Mehrzi, F. al-Saeedi,  A Deformed Wave Equation and Huygens’ Principle, Mathematics, 8(1), 10 (2019).

\bibitem{BD} S. Ben Sa\"id, L. Deleaval, Translation Operator and Maximal Function for the $(k,1)$-Generalized Fourier Transform, Journal of Functional Analysis, vol. 279, no. 8 (2020), 108706.

\bibitem{BH} 
W. Bloom and H. Heyer, Harmonic Analysis of Probability Measures on Hypergroups. De
Gruyter-Verlag, Berlin, 1994. 

\bibitem{BNS}M. A. Boubatra, S. Negzaoui, M. Sifi, A new product formula involving Bessel functions, Integral Transforms and Special Functions, (2022).

\bibitem{BO} S. Ben Sa\"id, B. \O rsted, The wave equation for Dunkl operators, Indagationes Mathematicae 16.3-4 (2005), 351--391.

\bibitem{BSKO} S. Ben Sa\"id, T. Kobayashi and B. \O rsted, Laguerre semigroup and Dunkl operators, Compos. Math., 148(4)(2012), 1265--1336.


\bibitem{CD} D. Constales, H. De Bie and P. Lian, Explicit formulas for the Dunkl dihedral kernel
and the $(k,a)$-generalized Fourier kernel, J. Math. Anal. Appl. 460(2) (2018), 900--926.

\bibitem{DH}
J. J. Duistermaat, L. Hörmander, Fourier Integral Operators. II. In: Brüning, J., Guillemin, V.W. (eds) Mathematics Past and Present Fourier Integral Operators. Springer, Berlin, Heidelberg (1994).


\bibitem{Du} C.F. Dunkl, Reflection groups and orthogonal polynomials on the
sphere, Math. Z., 197
 (1988), 33--60.

\bibitem{Du1} C.F. Dunkl, Differential-difference operators associated to reflection groups, Trans. Amer. Math. Soc.,  311, no. 1(1989), 167--183.

\bibitem{Du2} C.F. Dunkl, Hankel transforms associated to finite reflection groups, Hypergeometric Functions on Domains of Positivity, Jack Polynomials, and Applications,  Proceedings of an AMS Special Session Held March 22-23, 1991 in Tampa, Florida, Vol. 138, p. 123, American Mathematical Soc..


\bibitem{DW} F. Dai, H. Wang, A transference theorem for the Dunkl transform and its applications, J. Funct.
Anal. 258.12(2010), 4052--4074.

\bibitem{Er}
A. Erd\'elyi, et al, Tables of Integral Transforms, vol. II, McGraw Hill, New York, 1954. 

\bibitem{FN}
N. Faustino, S. Negzaoui, A Novel Schwartz Space for the $\left(k,\frac{2}{n}\right)-$Generalized Fourier Transform, arXiv:2507.04064.


\bibitem{GIT2} D. Gorbachev, V. Ivanov, S. Tikhonov, Pitt's inequalities and uncertainty principle for generalized Fourier transform, Int. Math. Res. Not., Issue 23(2016), 7179--7200.


\bibitem{GIT3} D.V. Gorbachev, V.I. Ivanov, S.Yu. Tikhonov, On the kernel of the $(\kappa,a)$-generalized Fourier transform, 	Forum of Mathematics, Sigma, 2023, 11:e72.. 

\bibitem{Ho} L. H\"ormander, The Analysis of Linear Partial Differential Operators I. Springer-Verlag,
Berlin, 1983.

\bibitem{H} R. Howe, The oscillator semigroup. The mathematical heritage of Hermann Weyl (Durham, NC, 1987), 61--132, Proc. Sympos. Pure Math., 48, Amer. Math. Soc., Providence, RI, 1988.

\bibitem{Iv}
V.I. Ivanov, One-dimensional $(k,a)$-generalized Fourier transform, Trudy Instituta Matematiki i Mekhaniki UrO RAN, 2023, vol. 29, no. 4, pp. 92--108.

\bibitem{Iv1}
 V.I.Ivanov, Bounded translation operator for the $(k, 1)$-generalized Fourier transform. Chebyshevskii Sbornik. 2020;21(4):85-96. (In Russ.)



\bibitem{Jo} T. R. Johansen, Weighted inequalities and uncertainty principles for the $(k,a)$-generalized Fourier transform, Int. J. Math.
27(3) (2016), 1650019.


\bibitem{KM1} T. Kobayashi, G. Mano, The inversion formula and holomorphic extension of the minimal representation of the conformal group, Harmonic Analysis, Group Representations, Automorphic Forms and
Invariant Theory: In honor of Roger Howe, (eds. J. S. Li, E. C.
Tan, N. Wallach and C. B. Zhu), World Scientific(2007), 159--223.

\bibitem{KM2} T. Kobayashi, G. Mano, The Schr\"odinger model for the minimal representation of the indefinite
orthogonal group $O(p,q)$, Mem. Amer. Math. Soc. 213(1000), 2011.

\bibitem{Me1}
H. Mejjaoli, K. Trim\`eche, On a mean value property associated with the Dunkl Laplacian
operator and applications, Integral Transforms and Special Functions, 12:3 (2001), 279--302.


\bibitem{Me}
H.Mejjaoli, Generalized translation operator and uncertainty principles associated with the deformed Stockwell transform, Revista de la Uni\'on Matemática Argentina, 2023, 65(2): 375--423.


\bibitem{Mu2} B. Muckenhoupt, Mean convergence of Hermite and Laguerre series II, ibid. 147 (1970),
433--460.

\bibitem{No} A. Nowak, Heat-diffusion and Poisson integrals for Laguerre and special Hermite expansions on weighted $L^p$ spaces. Studia Mathematica, 158(3)(2003), 239--268.


\bibitem{R1} M. R\"osler, Bessel-type signed hypergroups on $\mathbb R$, in: Probability Measures on Groups and Related Structures, XI, Oberwolfach, 1994, World Sci. Publ., River Edge, NJ, 1995.



\bibitem{R3} M. R\"osler, Positivity of Dunkl's intertwining operator. Duke Mathematical Journal, vol. 98, no. 3 (1999) pp. 445--463.


\bibitem{R2} M. R\"osler, A positive radial product formula for the Dunkl kernel, Trans. Amer.Math. Soc. 355(2003), no. 6, 2413--2438.

\bibitem{Te2} W. Teng, Hardy Inequalities for Fractional $(k, a)$-Generalized Harmonic Oscillators, Journal of Lie Theory, 32(4)(2022), 1007--1023.

\bibitem{Te3}
W. Teng, Imaginary powers of $(k,1)$-generalized harmonic oscillator, Complex Anal. Oper. Theory 16, 89 (2022).

\bibitem{T} K. Trim\`eche, Transformation int\'egrale de Weyl
   et th\'eor\`eme de Paley-Wiener associ\'es \`a un op\'erateur
   diff\'erentiel singulier sur $(0,\infty)$,
J. Math. Pures et Appl, 60 (1981), 51--98.



\end{thebibliography}
\end{document}